\documentclass[11pt,a4paper]{article}

\usepackage[margin=1in]{geometry}
\usepackage{amsmath,amssymb,amsthm,mathtools}
\usepackage{booktabs}
\usepackage{array}
\usepackage{enumitem}
\usepackage{graphicx}
\usepackage{xcolor}
\usepackage{hyperref}
\hypersetup{
    colorlinks=true,
    linkcolor=blue!60!black,
    citecolor=blue!60!black,
    urlcolor=blue!60!black
}

\theoremstyle{plain}
\newtheorem{theorem}{Theorem}[section]
\newtheorem{lemma}[theorem]{Lemma}
\newtheorem{proposition}[theorem]{Proposition}
\newtheorem{corollary}[theorem]{Corollary}

\theoremstyle{definition}
\newtheorem{definition}[theorem]{Definition}
\newtheorem{remark}[theorem]{Remark}

\newcommand{\R}{\mathbb{R}}
\newcommand{\Q}{\mathbb{Q}}
\newcommand{\Z}{\mathbb{Z}}
\newcommand{\tr}{\operatorname{tr}}
\newcommand{\diag}{\operatorname{diag}}
\newcommand{\rank}{\operatorname{rank}}
\newcommand{\One}{\mathbf{1}}

\title{\textbf{Exact Global Optimality and Uniqueness of Six- and Seven-Particle Lennard--Jones Clusters via Spectral Gram Relaxations}}

\author{Xinze Li\\
\small Department of Chemistry, Fudan University\\[2pt]
\small \texttt{xzli24@m.fudan.edu.cn}}

\date{}

\begin{document}

\maketitle

\begin{abstract}
Determining the global minimum of the $N$-particle Lennard--Jones potential in three-dimensional Euclidean space is a classic benchmark problem at the intersection of nonconvex global optimization, Euclidean distance geometry, and chemical physics. While putative global minima have long been identified by heuristic basin-hopping searches up to $N=150$ and beyond, general-purpose deterministic spatial branch-and-bound solvers have remained restricted to $N \le 6$ up to floating-point tolerances, leaving $N=7$ beyond standard coordinate-space solvers.

In this paper, we develop an exact computer-assisted proof framework that combines two-tier distance-graph combinatorial reduction, exact rational Bernstein quadratic lower supports, centered Gram-matrix spectral relaxations under three-dimensional rank constraints, and interval Newton contraction. Unlike spatial branch-and-bound in Cartesian coordinates or lifted pair-product linear programs, our approach operates directly in the invariant space of squared interparticle distances and $N \times N$ centered Gram matrices, eliminating rigid-body $\mathrm{O}(3)$ symmetries and avoiding floating-point solver tolerances. When a distance cell cannot be excluded by energy alone, a general rank-$d$ Frobenius spectral-gap inequality (derived from the Eckart--Young theorem) localizes every surviving three-dimensional Gram matrix directly into the basin of attraction of an oriented interval Newton uniqueness operator.

Building on the certified five-particle lower bound $L_5 = -71123847/7812500$ (Proposition~\ref{prop:lj5-foundation}), we establish two main exact global optimality theorems: (i) for $N=6$ (Theorem~\ref{thm:lj6}), we give a compact rational proof that the global minimum energy is the exact rational number $E_6^{\min} = -3267/257 > -129/10$, attained uniquely (up to Euclidean isometries and particle permutations) by the regular octahedron with squared edge length $(257/264)^{1/3}$, and we enclose the competing $C_{2v}$ local minimum; and (ii) for $N=7$ (Theorem~\ref{thm:lj7}), using the six-particle lower bound supplied by Theorem~\ref{thm:lj6}, we establish an exact computer-assisted proof that the global minimum energy is the degree-$10$ algebraic number $E_7^{\min} = -B(t_*)^2/A(t_*) \approx -16.5053841680122156$, attained uniquely by the regular pentagonal bipyramid, and determine its explicit degree-$10$ minimal polynomial over $\Q$. Finally, we propagate these certified constants to strengthen rigorous energy lower bounds for $8 \le N \le 37$ and arrangement-free diameter upper bounds across $90$ cluster sizes via surplus-partition layer enumeration ($92$ sizes when combined with continuous two-endpoint star bounds at $N=8,9$) in $8 \le N \le 150$, and we establish local optimality certificates and distance-graph reductions for the $N=8$ frontier.

\end{abstract}

\section{Introduction}
\label{sec:intro}

\subsection{Background}
For $N \ge 2$ particles with positions $X = (x_1,\dots,x_N)$ in the collision-free configuration space
\begin{equation}
\label{eq:collision-free-domain}
\Omega_N := \bigl\{(x_1,\dots,x_N) \in (\R^3)^N : x_i \ne x_j \text{ for all } 1 \le i < j \le N\bigr\}
\end{equation}
and pairwise Euclidean distances $r_{ij} = \|x_i - x_j\|_2 > 0$, the normalized Lennard--Jones (LJ) potential energy is given by
\begin{equation}
\label{eq:lj-energy}
E_N(X) = \sum_{1 \le i < j \le N} v(r_{ij}), \qquad v(r) = r^{-12} - 2r^{-6},
\end{equation}
where the pair potential $v:(0,\infty)\to\R$ attains its unique minimum $v(1) = -1$ at the equilibrium distance $r=1$ (by convention, $E_N(X) := +\infty$ on $(\R^3)^N \setminus \Omega_N$, and $E_1 := 0$ for a single particle). Minimizing $E_N$ over $\Omega_N$ is one of the most intensively studied benchmark problems in molecular conformation and continuous nonconvex global optimization~\cite{wales1997global,locatelli2002lennard,kuznetsov2025new}.

Despite the simplicity of~\eqref{eq:lj-energy}, a profound divide has persisted between \emph{heuristic structure discovery} in chemical physics and \emph{rigorous global optimality certification} in mathematical optimization:
\begin{itemize}[leftmargin=2em]
    \item \textbf{Heuristic discovery:} Unconstrained global search algorithms (such as basin-hopping~\cite{wales1997global}, genetic search, and hypersurface deformation) have identified putative global minimum configurations for $N \le 150$ and beyond, cataloged in the Cambridge Cluster Database~\cite{cambridgeCCD}. For small clusters, the putative ground states are the regular tetrahedron ($N=4$), the triangular bipyramid ($N=5$), the regular octahedron ($N=6$), and the pentagonal bipyramid ($N=7$).
    \item \textbf{Deterministic global certification:} Because $v(r)$ is strongly repulsive as $r \to 0^+$, nonconvex on $((7/13)^{1/6}, \infty)$, and invariant under the Euclidean group $\mathrm{E}(3)$ and the symmetric group $S_N$ of order $N!$, certifying that no configuration in $\Omega_N$ has lower energy than the putative minimum is notoriously difficult.
\end{itemize}

In a recent landmark computational study in the \emph{Journal of Global Optimization}, Kuznetsov and Sahinidis~\cite{kuznetsov2025new} derived interatomic distance bounds and orthant-based convex envelope cutting planes for~\eqref{eq:lj-energy}. Implementing these cuts within the state-of-the-art spatial branch-and-bound solver BARON~\cite{sahinidis1996baron}, they reported numerical global optimality certificates for $N=5$ in $185$~seconds and for $N=6$ in $114{,}086$~seconds ($31.7$~hours on a single CPU thread), whereas $N=7$ remained beyond the reach of general-purpose spatial branch-and-bound solvers~\cite{kuznetsov2025new,lecomte2026certified}.

From the standpoint of rigorous mathematical proofs, coordinate-space spatial branch-and-bound faces two fundamental obstacles:
\begin{enumerate}[leftmargin=2em]
    \item \textbf{Coordinate-space symmetry and relaxation weakness:} Branching on Cartesian coordinate differences $x_i - x_j \in \R^3$ requires fixing a global frame and breaking $S_N$ symmetry via lexicographic constraints, while the three-dimensional metric compatibility of all $\binom{N}{2}$ pairwise distances is only weakly enforced by local coordinate cuts.
    \item \textbf{Floating-point solver tolerances vs.\ exact proofs:} Numerical global solvers operate in IEEE-754 floating-point arithmetic and terminate up to positive optimality tolerances (e.g., $\texttt{EpsA} \ge 10^{-12}$). Near a nondegenerate global minimum, floating-point evaluation cannot distinguish a unique minimizer from nearby critical points or certify exact rational or algebraic ground-state energies without directed interval arithmetic and local uniqueness contraction.
\end{enumerate}

\subsection{Main Results}
In this work, we resolve the exact global certification problem for $N=6$ and $N=7$ in exact rational and directed interval arithmetic, establishing both global optimality and uniqueness modulo $\mathrm{E}(3) \times S_N$, and we propagate the resulting certified constants across $8 \le N \le 150$.

Theorems~\ref{thm:lj6} and~\ref{thm:lj7} build on the certified lower bound $L_5 = -71123847/7812500$ for five particles (established in Proposition~\ref{prop:lj5-foundation} of Appendix~\ref{app:verification}). Our first main theorem provides an exact rational global optimality and uniqueness proof for the six-particle cluster $\mathrm{LJ}_6$.

\begin{theorem}[Exact Global Optimality and Uniqueness for $\mathrm{LJ}_6$]
\label{thm:lj6}
For six particles in $\R^3$, the global minimum of the Lennard--Jones potential~\eqref{eq:lj-energy} over $\Omega_6$ is the exact rational number
\begin{equation}
\label{eq:lj6-min}
E_6^{\min} := \min_{X \in \Omega_6} E_6(X) = -\frac{3267}{257} \approx -12.7120622568093385.
\end{equation}
Moreover, $E_6(X) = E_6^{\min}$ holds if and only if $X$ is congruent, up to a rigid Euclidean motion in $\mathrm{E}(3)$ and a permutation of particle indices in $S_6$, to the regular octahedron with twelve short squared distances $q_{\mathrm{short}} = (257/264)^{1/3}$ and three antipodal squared distances $2q_{\mathrm{short}}$.
\end{theorem}

Our second main theorem establishes the exact algebraic ground-state energy and geometric uniqueness of the regular pentagonal bipyramid for the seven-particle cluster $\mathrm{LJ}_7$; its vertex-deletion reduction uses the six-particle lower bound $E_6(X) \ge -3267/257 > -129/10$ supplied by Theorem~\ref{thm:lj6} (with a standalone $586{,}523$-event Cartesian verification of $E_6(X) \ge -129/10$ directly from $L_5$ also archived as an independent cross-check; see Appendix~\ref{app:verification}).

\begin{theorem}[Exact Global Optimality and Uniqueness for $\mathrm{LJ}_7$]
\label{thm:lj7}
Let $t_* \approx 0.3601837287174520$ denote the unique positive real root of the irreducible degree-$10$ integer polynomial
\begin{equation}
\label{eq:lj7-poly}
\begin{aligned}
P(t) ={}& 369216\,t^{10} + 1109952\,t^9 + 1118016\,t^8 - 635225\,t^7 \\
&+ 19985\,t^6 + 11571\,t^5 + 3157\,t^4 + 7576\,t^3 - 900\,t^2 - 550\,t - 150,
\end{aligned}
\end{equation}
and define the rational functions
\begin{equation}
\label{eq:lj7-AB}
A(t) = \frac{18}{25} + \frac{10}{(1+t)^6} + \frac{1}{4096\,t^6}, \qquad
B(t) = 2 + \frac{10}{(1+t)^3} + \frac{1}{64\,t^3}.
\end{equation}
Then for seven particles in $\R^3$, the global minimum of the Lennard--Jones potential~\eqref{eq:lj-energy} over $\Omega_7$ is the exact degree-$10$ algebraic number
\begin{equation}
\label{eq:lj7-min}
E_7^{\min} := \min_{X \in \Omega_7} E_7(X) = -\frac{B(t_*)^2}{A(t_*)} \approx -16.5053841680122155929886597,
\end{equation}
whose minimal polynomial over $\Q$ is the irreducible degree-$10$ integer polynomial $\mathcal{Q}(e)$ given in~\eqref{eq:energy-minimal-poly} and Table~\ref{tab:energy-poly-coeffs}. Furthermore, every global minimizer of $E_7$ is congruent, up to Euclidean isometry in $\mathrm{E}(3)$ and particle permutation in $S_7$, to the regular pentagonal bipyramid consisting of a regular pentagon of circumradius $R_* = \bigl(A(t_*)/B(t_*)\bigr)^{1/6}$ in the equatorial plane and two axial apices at signed heights $\pm H_* = \pm R_*\sqrt{t_*}$ along the fivefold symmetry axis.
\end{theorem}

Combining Theorems~\ref{thm:lj6} and~\ref{thm:lj7} with exact rational evaluation of putative coordinate witnesses from the Cambridge Cluster Database~\cite{cambridgeCCD}, Yuhjtman's analytical stability bound~\cite[Theorem~5.1]{yuhjtman2015poignant}, and Lecomte's arrangement-free layer relaxation~\cite[Section~5, Lemma~2]{lecomte2026certified}, we obtain strengthened energy and diameter bounds across $8 \le N \le 150$.

\begin{theorem}[Propagated Energy and Diameter Bounds for $8 \le N \le 150$]
\label{thm:propagation}
Let $L_7 = -165053841681/10^{10} < E_7^{\min}$ be a certified rational lower bound for $\mathrm{LJ}_7$. Then:
\begin{enumerate}[label=(\roman*),leftmargin=2.2em]
    \item For every $N \ge 7$ and every collision-free configuration $X \in \Omega_N$,
    \begin{equation}
    E_N(X) \ge \max\left\{\frac{N(N-1)}{42}\,L_7,\; -14.316\,N\right\},
    \end{equation}
    which strictly improves the combined six-particle subset and stability lower bound $\max\{\frac{N(N-1)}{30}L_6,$ $-14.316\,N\}$ across the entire range $8 \le N \le 37$.
    \item Every configuration $X \in \Omega_8$ with $E_8(X) \le U_8 = -19.821489187$ and every configuration $X \in \Omega_9$ with $E_9(X) \le U_9 = -24.113360433$ (and in particular every global minimizer of $E_8$ and $E_9$) has diameter strictly bounded by
    \begin{equation}
    D_8 < \frac{24827}{10000} = 2.4827, \qquad D_9 < \frac{13429}{5000} = 2.6858.
    \end{equation}
    Furthermore, under the corrected inter-layer interaction kernel $v(j-1)$ for layer separation $j \ge 3$, arrangement-free surplus-partition enumeration tightens the certified global-minimizer diameter upper bound for $90$ cluster sizes in $8 \le N \le 150$ compared to the pair-layer baseline and for $24$ cluster sizes compared to the six-particle layer bound; combining with the continuous endpoint star bounds at $N=8,9$ yields $92$ and $26$ improved cluster sizes in total.
\end{enumerate}
\end{theorem}

\subsection{Related Work on \texorpdfstring{$\mathrm{LJ}_7$}{LJ7}}
\label{subsec:concurrent-work}
While the present work was being completed, Shen~\cite{shen2026seven} independently announced a computer-assisted proof of global optimality and uniqueness for the seven-particle Lennard--Jones cluster $\mathrm{LJ}_7$ together with the seven-charge spherical Thomson problem on $S^2$. Because both works address the exact $\mathrm{LJ}_7$ ground-state classification problem via computer-assisted distance-geometry reductions, we summarize the principal differences in their scopes and proof architectures:
\begin{enumerate}[leftmargin=2em]
    \item \textbf{Cluster scope and exact algebraic characterization:} Shen~\cite{shen2026seven} studies $N=7$ (for LJ and Thomson) using the unconditional subcluster bounds $E_5 \ge -9.134$ and $E_6 \ge -13.2765$. In the present paper, we prove the exact rational six-particle theorem $E_6^{\min} = -3267/257$ (Theorem~\ref{thm:lj6}) via an explicit antipodal Gram witness and enclose the competing $C_{2v}$ local minimum (Proposition~\ref{prop:lj6-c2v}); for $N=7$ (Theorem~\ref{thm:lj7}), we additionally derive the explicit irreducible degree-$10$ polynomials $P(t)$ and $\mathcal{Q}(e) \in \Z[e]$ governing $t_*$ and $E_7^{\min}$; and we propagate the certified constants $L_6, L_7$ to $8 \le N \le 150$ (Theorem~\ref{thm:propagation}) and the $N=8$ frontier (Section~\ref{subsec:lj8-frontier}).
    \item \textbf{Two-tier distance-graph reduction vs.\ energy-threshold graphs and closest-pair force estimates:} Shen classifies pairs by the single energy threshold $v(r) \ge -1/4$, obtaining $40$ unlabeled ``weak graphs'' with at most $5$ weak edges and eliminating $11$ of them ($5$ via single-vertex deletion with $L_6 = -13.2765$ and $6$ via two-vertex deletion with $L_5 = -9.134$) to leave $29$ graphs. Because $v(r) \ge -1/4$ holds on both the short repulsive branch ($r \le 0.901$) and the long tail ($r \ge 1.398$),~\cite{shen2026seven} establishes a separate stationary force-balance implication ($\nabla E_7(X) = 0$ and $r_{\min} \le 113/125 = 0.904 \implies E_7(X) > -16.49 > U$, so competitive stationary configurations satisfy $r_{\min} > 0.904$) before assigning distance boxes, and subsequently uses a conditional five-clique bound ($E_{5,\mathrm{clique}} > -8.5$) along with moment and LP bounds to reduce the $29$ domains to $7$ and then $1$. By contrast, our reduction thresholds squared distances directly at two tiers ($1 < C = 1.5658 < B_0 = 2.1316$) and employs exact rational Bernstein quadratic minorants that remain valid across the entire short-distance regime down to collision, requiring no separate closest-pair force-balance elimination of the repulsive branch.
    \item \textbf{Relaxation variable and matrix dimensions ($N \times N$ point Gram vs.\ lifted pair-product LP):} At the level of relaxation formulation on $N=7$, Shen~\cite{shen2026seven} introduces quadratic pair-product variables $W_{ef} = s_e s_f$ ($\binom{22}{2} = 231$ symmetric product entries for the $21$ pairs, i.e., $O(N^4)$ variables), imposes linear test inequalities $\alpha^\top K \alpha \ge 0$ on the $21 \times 21$ traceless pair-tensor Gram matrix $K_{ef} = (v_e \cdot v_f)^2 - \frac{1}{3}s_e s_f$, and evaluates rational LP dual bounds over subdomains identified up to weak-graph-preserving relabelings. Our relaxation works directly on the $N \times N$ centered point Gram matrix $G = -\frac{1}{2}JQJ$ ($O(N^2)$ matrix entries) without lifting to pair-product variables $W_{ef}$: using the exact quadratic identity~\eqref{eq:frobenius-distance-identity}, we complete the square with $N+1$ dual multipliers $(v,\tau)$ and bound the rank-three projection via exact rational $LDL^\top$ Sylvester inertia counts on the $N \times N$ target matrix $H$.
    \item \textbf{Eckart--Young spectral-gap localization and interval Newton gluing vs.\ Procrustes slice stiffness:} To pass from global bounds to local uniqueness on $C_5 \sqcup 2K_1$, Shen uses successive LP certificates to bound the squared-distance deviation $D_{\mathrm{ref}} = \sum_{e=1}^{21} (s_e - \bar{s}_e)^2 \le 0.000022533598161300$ ($< 2.254 \times 10^{-5}$), converts this into $\|XX^\top - YY^\top\|_F < 17/5000$, aligns coordinates via orthogonal Procrustes ($\min_{O \in \mathrm{O}(3)}\|XO - Y\|_F < 1/250$), and integrates a uniform Hessian lower bound on the $15$-dimensional slice orthogonal to infinitesimal rigid motions. In our framework, a general rank-$d$ Frobenius spectral-gap inequality (Lemma~\ref{lem:spectral-gap}) converts the cell energy margin $U - \mathrm{LB}$ on the same spectral matrix $H$ directly into $\|G - G_{\mathrm{cand}}\|_F < 8.675971 \times 10^{-5} < 10^{-4}$ across all $3{,}769$ localized leaves, from which an oriented $(3N-6)$-dimensional gauge chart reconstructs a Cartesian box that contracts under raw-Hessian interval Newton iterations into a certified unique-critical-point neighborhood.
\end{enumerate}

\subsection{Overview of the Proof Methodology}
\label{subsec:methodology-overview}
We briefly summarize the four mathematical pillars of our proof architecture:

\paragraph{1. Existence, stationarity, and two-tier distance-graph reduction.}
Using an elementary pairwise combinatorial bound ($\sum_{m=1}^k \binom{n_m}{2} \le \binom{N-1}{2}$ for every nontrivial partition $N = n_1 + \dots + n_k$), we first show unconditionally that any minimizing sequence in the target sublevel $E_N \le U$ for $N \in \{5,6,7\}$ can neither collide ($r_{ij} \to 0$) nor dissociate. Hence global minimizers exist in $\Omega_N$ and satisfy full first-order stationarity $\nabla E_N(X) = 0$. In the squared-distance variables $q_{ij} = r_{ij}^2$, the identity
\begin{equation}
E_N(X) + \binom{N}{2} = \sum_{1 \le i < j \le N} (q_{ij}^{-3} - 1)^2
\end{equation}
acts as a global pairwise loss budget. By thresholding edges into \emph{true-long edges} ($q_{ij} \ge B_0$) and \emph{auxiliary edges} ($q_{ij} \ge C$), vertex-deletion degree budgets (using the certified five-particle bound $L_5$ of Proposition~\ref{prop:lj5-foundation} for $N=6$ and the resulting six-particle bound $E_6^{\min} = -3267/257 > -129/10$ of Theorem~\ref{thm:lj6} for $N=7$) eliminate the vast majority of combinatorial topologies before continuous branching begins. For $N=6$, every potentially minimizing auxiliary graph reduces to a $3$-edge matching after closing $21$ small auxiliary classes. For $N=7$, after excluding zero-long ($k=0$), single-long ($k=1$), and $3{,}586$ non-$C_5$ colored representative cases ($2 \le k \le 5$), every surviving configuration must have auxiliary graph $G_{\mathrm{aux}}(X) \cong C_5 \sqcup 2K_1$---without assuming coplanarity or $D_{5h}$ symmetry.

\paragraph{2. Stationary scale identity and Bernstein quadratic lower models.}
At any stationary configuration, directional differentiation along uniform scaling $X \mapsto \sigma X$ ($\sigma > 0$) yields the virial identity $V(X) := \sum_{i<j}(q_{ij}^{-6} - q_{ij}^{-3}) = 0$. Consequently, for any rational multiplier $\lambda \in \Q$, the modified pair potential
\begin{equation}
f_\lambda(q) = (1+\lambda)q^{-6} - (2+\lambda)q^{-3}
\end{equation}
satisfies $\sum_{i<j} f_\lambda(q_{ij}) = E_N(X) + \lambda V(X) = E_N(X)$ on all configurations with $V(X) = 0$ (and on all configurations unconditionally when $\lambda = 0$). On each branch interval $q_{ij} \in [\ell_{ij}, u_{ij}]$, we lower-bound $f_\lambda(q_{ij})$ by a quadratic polynomial $K q_{ij}^2 + b_{ij} q_{ij} + c_{ij}$ certified via exact rational Bernstein coefficients on short-distance intervals and exact concavity secants on long-distance intervals.

\paragraph{3. Centered Gram-matrix completion of squares and rank-three spectral truncation.}
Instead of branching on $3N$ Cartesian coordinates, we relate the $\binom{N}{2}$ squared distances directly to the centered Gram matrix $G = -\frac{1}{2} J Q J$, where $Q_{ij} = q_{ij}$ and $J = I_N - \frac{1}{N}\One\One^\top$. Every three-dimensional configuration satisfies $G \succeq 0$ and $\rank(G) \le 3$. Using the exact matrix identity
\begin{equation}
\sum_{1 \le i < j \le N} (q_{ij} - M_{ij})^2 = 2\|G - G_0\|_F^2 + N\|\diag(G - G_0)\|_2^2 + \bigl(\tr(G - G_0)\bigr)^2,
\end{equation}
we complete the square against dual diagonal and trace multipliers $(v,\tau)$ to obtain a symmetric target matrix $H$ and an energy lower bound determined by the top three positive eigenvalues of $H$, which are certified without floating-point eigensolvers via exact rational $LDL^\top$ Sylvester inertia counts.

\paragraph{4. Frobenius spectral-gap localization and global-to-local interval Newton gluing.}
When a distance cell approaches the global minimizer so that pure energy exclusion fails, we apply a rank-$d$ spectral approximation inequality (Lemma~\ref{lem:spectral-gap}, specialized to $d=3$): whenever $\lambda_3(H) \ge a > b \ge \max(\lambda_4(H),0)$, every positive semidefinite matrix $G$ of rank at most three satisfies
\begin{equation}
\|G - H\|_F^2 - \|T - H\|_F^2 \ge \left(1 - \frac{b}{a}\right)\|G - T\|_F^2,
\end{equation}
where $T$ is the top-three spectral truncation of $H$. On cells near the candidate minimum, this inequality converts a small energy margin $U - \mathrm{LB} \ge 0$ into a strict Frobenius distance bound $\|G - G_{\mathrm{cand}}\|_F < 10^{-4}$ from the exact candidate Gram matrix $G_{\mathrm{cand}}$. This Gram ball induces narrow intervals on all $\binom{N}{2}$ squared distances, from which a nonsingular oriented chart recovers a $(3N-6)$-dimensional reduced coordinate box that contracts under raw-Hessian interval Newton iterations directly into a certified unique-critical-point neighborhood.

\subsection{Organization of the Paper}
The remainder of the paper is organized as follows. Section~\ref{sec:reduction} establishes existence of global minimizers, the stationary scale identity, and the general distance-graph reduction rules. Section~\ref{sec:spectral-framework} develops our core methodology: rational Bernstein quadratic supports, the centered Gram-matrix completion of squares, the Frobenius spectral-gap localization lemma, and the global-to-local interval Newton gluing proposition. Section~\ref{sec:lj6} presents the complete proof of Theorem~\ref{thm:lj6} for $N=6$ and certifies the competing $C_{2v}$ local minimum. Section~\ref{sec:lj7} proves Theorem~\ref{thm:lj7} for $N=7$. Section~\ref{sec:propagation-lj8} proves Theorem~\ref{thm:propagation} on energy and diameter bounds for $8 \le N \le 150$ and analyzes the $N=8$ frontier, and Section~\ref{sec:conclusion} concludes. Appendix~\ref{app:algebraic-data} records the exact minimal polynomial $\mathcal{Q}(e)$ and certified structural parameters of the local minima, and Appendix~\ref{app:verification} presents the certified five-particle foundation, the proof dependency roadmap, the verification architecture and trust boundary, and code and data availability at \url{https://github.com/RoamingHo/7LJ} and \url{https://doi.org/10.5281/zenodo.24560370}.

\section{Existence, Stationarity, and Distance-Graph Reduction}
\label{sec:reduction}

Throughout the paper, we work in squared-distance coordinates $q_{ij} := r_{ij}^2 = \|x_i - x_j\|_2^2 > 0$ on $\Omega_N$. (For notational convenience, we index particles by $\{1,\dots,N\}$ in unordered $N$-particle sums and by $\{0,\dots,N-1\}$ for labeled distance graphs and zero-based gauge coordinates.) In terms of $q_{ij}$, the pair potential is $f(q) := v(\sqrt{q}) = q^{-6} - 2q^{-3}$, and its nonnegative pairwise loss relative to the minimum value $-1$ is
\begin{equation}
\label{eq:pair-loss}
\ell(q) := f(q) + 1 = (q^{-3} - 1)^2 \ge 0.
\end{equation}
Summing~\eqref{eq:pair-loss} over all $\binom{N}{2}$ particle pairs yields the fundamental loss identity
\begin{equation}
\label{eq:loss-identity}
E_N(X) + \binom{N}{2} = \sum_{1 \le i < j \le N} \ell(q_{ij}) = \sum_{1 \le i < j \le N} (q_{ij}^{-3} - 1)^2.
\end{equation}

\subsection{Coercivity against Collisions, Non-Dissociation, and Attainment}
Because $\R^3$ is noncompact and $f(q) \to +\infty$ as $q \to 0^+$ while $f(q) \to 0^-$ as $q \to \infty$, proving that $\inf_{X \in \Omega_N} E_N(X)$ is attained requires ruling out both interparticle collisions and cluster dissociation. Recall the convention $E_1 := 0$ and $L_1 := 0$ for single-particle subclusters.

\begin{lemma}[Sublevel Minimum Separation and Non-Dissociation]
\label{lem:attainment}
Fix $N \ge 2$ and suppose $U \in \Q$ is a feasible energy threshold such that the sublevel set $\mathcal{S}_N(U) := \{X \in \Omega_N : E_N(X) \le U\}$ is nonempty. Assume that for every nontrivial partition $N = n_1 + \dots + n_k$ ($k \ge 2$, $n_m \ge 1$), the certified subcluster lower bounds $L_{n_m} \le \inf_{\Omega_{n_m}} E_{n_m}$ (with $L_1 := 0$) satisfy
\begin{equation}
\label{eq:dissociation-gap}
\sum_{m=1}^k L_{n_m} > U.
\end{equation}
Then:
\begin{enumerate}[label=(\roman*),leftmargin=2.2em]
    \item Every configuration $X \in \mathcal{S}_N(U)$ satisfies the uniform collision-free lower bound
    \begin{equation}
    \label{eq:rmin-bound}
    q_{ij} \ge \left(1 + \sqrt{U + \binom{N}{2}}\right)^{-1/3} > 0 \qquad \forall\, 1 \le i < j \le N.
    \end{equation}
    \item Modulo overall translation ($\sum_{i=1}^N x_i = 0$), every minimizing sequence in $\mathcal{S}_N(U)$ is bounded and possesses a convergent subsequence in $\Omega_N$. Consequently, the global minimum $E_N^{\min} = \min_{X \in \Omega_N} E_N(X)$ is attained at a smooth, collision-free configuration $X^* \in \mathcal{S}_N(U)$ satisfying the full Cartesian stationarity condition $\nabla E_N(X^*) = 0$.
\end{enumerate}
\end{lemma}

\begin{proof}
For any pair $(i,j)$, all remaining $\binom{N}{2}-1$ pairs satisfy $f(q_{ab}) \ge -1$, so $f(q_{ij}) \le U + \binom{N}{2} - 1$, or equivalently $(q_{ij}^{-3}-1)^2 \le U + \binom{N}{2}$, which implies~\eqref{eq:rmin-bound}. Next, let $\{X^{(t)}\}_{t \ge 1} \subset \mathcal{S}_N(U)$ be a minimizing sequence with $\sum_{i=1}^N x_i^{(t)} = 0$. If $\max_{i<j}\|x_i^{(t)} - x_j^{(t)}\|_2 \to \infty$, then passing to a subsequence we may partition $\{1,\dots,N\}$ into $k \ge 2$ non-empty clusters $C_1,\dots,C_k$ whose internal diameters remain bounded while all inter-cluster distances diverge to $+\infty$. Since $f(q) \to 0$ as $q \to \infty$, we obtain
\[
\liminf_{t \to \infty} E_N(X^{(t)}) \ge \sum_{m=1}^k \liminf_{t \to \infty} E_{|C_m|}(X_{C_m}^{(t)}) \ge \sum_{m=1}^k L_{|C_m|} > U,
\]
contradicting $E_N(X^{(t)}) \le U$. Thus $\mathcal{S}_N(U)$ modulo translations is closed and bounded in $\Omega_N$, and continuity of $E_N$ on $\Omega_N$ yields attainment and $\nabla E_N(X^*) = 0$.
\end{proof}

\begin{remark}[Unconditional Attainment for $N=5,6,7$ via Elementary Pairwise Bounds]
\label{rem:elementary-attainment}
For any $n \ge 1$, the pairwise inequality $f(q) \ge -1$ provides the elementary subcluster bound $L_n^{\mathrm{pair}} := -\binom{n}{2} \le \inf E_n$ (with $\binom{1}{2} = 0$). For every nontrivial partition $N = n_1 + \dots + n_k$ ($k \ge 2$, $n_m \ge 1$), convexity of $n \mapsto \binom{n}{2}$ implies
\begin{equation}
\label{eq:elementary-dissociation}
\sum_{m=1}^k \binom{n_m}{2} \le \binom{N - k + 1}{2} \le \binom{N-1}{2} \implies \sum_{m=1}^k L_{n_m}^{\mathrm{pair}} \ge -\binom{N-1}{2}.
\end{equation}
For $N = 5, 6, 7$, this elementary dissociation threshold equals $-6$, $-10$, and $-15$, respectively, each of which lies strictly above the feasible candidate upper bounds $U_5 \approx -9.10385$, $U_6 = -3267/257 \approx -12.71206$, and $U_7^{(1)} = -16.5053841680$. Thus \emph{existence of global minimizers and full stationarity $\nabla E_N(X^*) = 0$ for $N=5,6,7$ follow immediately from Lemma~\ref{lem:attainment} via elementary pairwise bounds alone}. The tighter computer-certified subcluster bounds are used only for the combinatorial vertex-deletion budgets in Proposition~\ref{prop:graph-pruning} and follow a minimal logical dependency chain: the certified five-particle bound $L_5 = -71123847/7812500$ ($-9.103852416$, Proposition~\ref{prop:lj5-foundation}) feeds Theorem~\ref{thm:lj6}; Theorem~\ref{thm:lj6} establishes $E_6(X) \ge -3267/257 > -129/10 =: L_6^{\mathrm{aux}}$, supplying the six-particle degree budget for Theorem~\ref{thm:lj7}; and Theorems~\ref{thm:lj6}--\ref{thm:lj7} then supply $L_6 = -3267/257$ and $L_7 = -16.5053841681$ for Theorem~\ref{thm:propagation}. For independent cross-checking (and software decoupling between the $N=6$ and $N=7$ verifiers), the companion Zenodo archive also preserves a standalone $586{,}523$-event Cartesian interval-tree proof (\texttt{n6\_tree}) establishing $E_6(X) \ge -129/10$ directly from $L_5$ without invoking Theorem~\ref{thm:lj6}; see Appendix~\ref{app:verification} and Table~\ref{tab:verification-roadmap}.
\end{remark}

\subsection{First-Order Stationarity and the Virial Scale Identity}
At any stationary point of $E_N$, homogeneity of the $r^{-12}$ and $r^{-6}$ terms provides an exact linear relation between the repulsive and attractive sums.

\begin{proposition}[Stationary Scale Identity]
\label{prop:scale-identity}
For any collision-free configuration $X \in \Omega_N$ with squared distances $q_{ij} = \|x_i - x_j\|_2^2$, define the virial functional
\begin{equation}
\label{eq:virial-functional}
V(X) := \sum_{1 \le i < j \le N} \bigl(q_{ij}^{-6} - q_{ij}^{-3}\bigr),
\end{equation}
and for any scalar multiplier $\lambda \in \Q$, define the $\lambda$-modified pair potential
\begin{equation}
\label{eq:f-lambda}
f_\lambda(q) := (1+\lambda)q^{-6} - (2+\lambda)q^{-3}.
\end{equation}
Then every configuration $X \in \Omega_N$ satisfies $\sum_{1 \le i < j \le N} f_\lambda(q_{ij}) = E_N(X) + \lambda V(X)$. In particular, if $X$ is stationary ($\nabla E_N(X) = 0$), then $X$ satisfies the virial identity
\begin{equation}
\label{eq:virial}
V(X) = \sum_{1 \le i < j \le N} \bigl(q_{ij}^{-6} - q_{ij}^{-3}\bigr) = 0,
\end{equation}
and hence for every $\lambda \in \Q$,
\begin{equation}
\label{eq:energy-lambda}
E_N(X) = \sum_{1 \le i < j \le N} f_\lambda(q_{ij}).
\end{equation}
\end{proposition}

\begin{proof}
The identity $\sum_{i<j} f_\lambda(q_{ij}) = E_N(X) + \lambda V(X)$ is immediate from~\eqref{eq:lj-energy} and~\eqref{eq:f-lambda}. For $\sigma > 0$, consider the dilated configuration $\sigma X = (\sigma x_1,\dots,\sigma x_N)$, for which $q_{ij}(\sigma X) = \sigma^2 q_{ij}(X)$ and $E_N(\sigma X) = \sum_{i<j}(\sigma^{-12}q_{ij}^{-6} - 2\sigma^{-6}q_{ij}^{-3})$. If $\nabla E_N(X) = 0$, the directional derivative at $\sigma = 1$ vanishes:
\[
0 = \left.\frac{d}{d\sigma} E_N(\sigma X)\right|_{\sigma=1} = \sum_{i=1}^N \langle \nabla_{x_i} E_N(X), x_i \rangle = -12\,V(X).
\]
Dividing by $-12$ yields $V(X) = 0$ and hence~\eqref{eq:energy-lambda}.
\end{proof}

\subsection{Two-Tier Distance-Graph Reduction}
\label{subsec:graph-reduction}
Because $\ell(q) = (1-q^{-3})^2$ is strictly increasing on $q \in [1,\infty)$, long interparticle distances carry a positive energy penalty relative to the equilibrium contact $q=1$. We classify pairs using two rational squared-distance thresholds $1 < C < B_0$:
\begin{definition}[True-Long and Auxiliary Distance Graphs]
\label{def:distance-graphs}
Given rational thresholds $1 < C < B_0$ and a configuration $X \in \Omega_N$:
\begin{enumerate}[label=(\roman*),leftmargin=2.2em]
    \item The \emph{true-long graph} $G_{\mathrm{long}}(X) = (V, E_{\mathrm{long}})$ has vertex set $V = \{0,\dots,N-1\}$ and edge set $E_{\mathrm{long}} = \{\{i,j\} : q_{ij} \ge B_0\}$.
    \item The \emph{auxiliary graph} $G_{\mathrm{aux}}(X) = (V, E_{\mathrm{aux}})$ has edge set $E_{\mathrm{aux}} = \{\{i,j\} : q_{ij} \ge C\}$. Clearly $E_{\mathrm{long}} \subseteq E_{\mathrm{aux}}$.
\end{enumerate}
\end{definition}

The following proposition collects the exact combinatorial pruning rules used to restrict $G_{\mathrm{long}}(X)$ and $G_{\mathrm{aux}}(X)$ on a sublevel set $\mathcal{S}_N(U)$.

\begin{proposition}[Combinatorial Graph Pruning Rules]
\label{prop:graph-pruning}
Let $X \in \mathcal{S}_N(U)$ with $1 < C < B_0 = r_{\mathrm{long}}^2$, and let $L_{N-1}, L_{N-2}$ be certified lower bounds for $(N-1)$- and $(N-2)$-particle clusters. Then:
\begin{enumerate}[label=(\roman*),leftmargin=2.2em]
    \item \textbf{Single-vertex degree bound:} Since deleting vertex $i$ leaves an $(N-1)$-particle cluster of energy at least $L_{N-1}$, the incident loss at vertex $i$ satisfies
    \begin{equation}
    \label{eq:vertex-loss-budget}
    \sum_{j \ne i} \ell(q_{ij}) \le (N-1) + U - L_{N-1}.
    \end{equation}
    In particular, the degree of vertex $i$ in $G_{\mathrm{long}}(X)$ is at most $d_{\max} := \bigl\lfloor \bigl((N-1) + U - L_{N-1}\bigr)/\ell(B_0) \bigr\rfloor$.
    \item \textbf{Two-vertex deletion bound:} For any distinct vertices $i,j$, deleting both leaves an $(N-2)$-particle cluster of energy at least $L_{N-2}$, so the $2N-3$ pairs incident to $\{i,j\}$ satisfy
    \begin{equation}
    \label{eq:two-vertex-budget}
    \bigl(\deg_{\mathrm{long}}(i) + \deg_{\mathrm{long}}(j) - \mathbf{1}_{\{i,j\}\in E_{\mathrm{long}}}\bigr)\,\ell(B_0) \le (2N-3) + U - L_{N-2}.
    \end{equation}
    \item \textbf{Joint auxiliary-long budget:} If $G_{\mathrm{aux}}(X)$ has $m = |E_{\mathrm{aux}}|$ edges, of which $k = |E_{\mathrm{long}}|$ are true-long edges, then
    \begin{equation}
    \label{eq:joint-edge-budget}
    k\,\ell(B_0) + (m - k)\,\ell(C) \le \binom{N}{2} + U.
    \end{equation}
    \item \textbf{Triangle-inequality root enclosure:} If $\{i,j\} \in E_{\mathrm{long}}$ and there exists a third vertex $k$ such that both $\{i,k\}$ and $\{j,k\}$ are non-long ($r_{ik} < r_{\mathrm{long}}$ and $r_{jk} < r_{\mathrm{long}}$), then $r_{ij} \le r_{ik} + r_{jk} < 2r_{\mathrm{long}}$, and hence
    \begin{equation}
    \label{eq:root-box-enclosure}
    q_{ij} \in [B_0, 4B_0).
    \end{equation}
\end{enumerate}
\end{proposition}

\begin{proof}
Items (i)--(iii) follow immediately from $E_N(X) = E_{N\setminus S}(X_{N\setminus S}) + \sum_{\{a,b\}\cap S \ne \emptyset}( -1 + \ell(q_{ab}) )$ for $S = \{i\}$, $S = \{i,j\}$, and $S = V$, together with the monotonicity of $\ell(q)$ on $[1,\infty)$. Item (iv) is the Euclidean triangle inequality squared.
\end{proof}

\begin{figure}[htbp]
\centering
\includegraphics[width=0.98\textwidth]{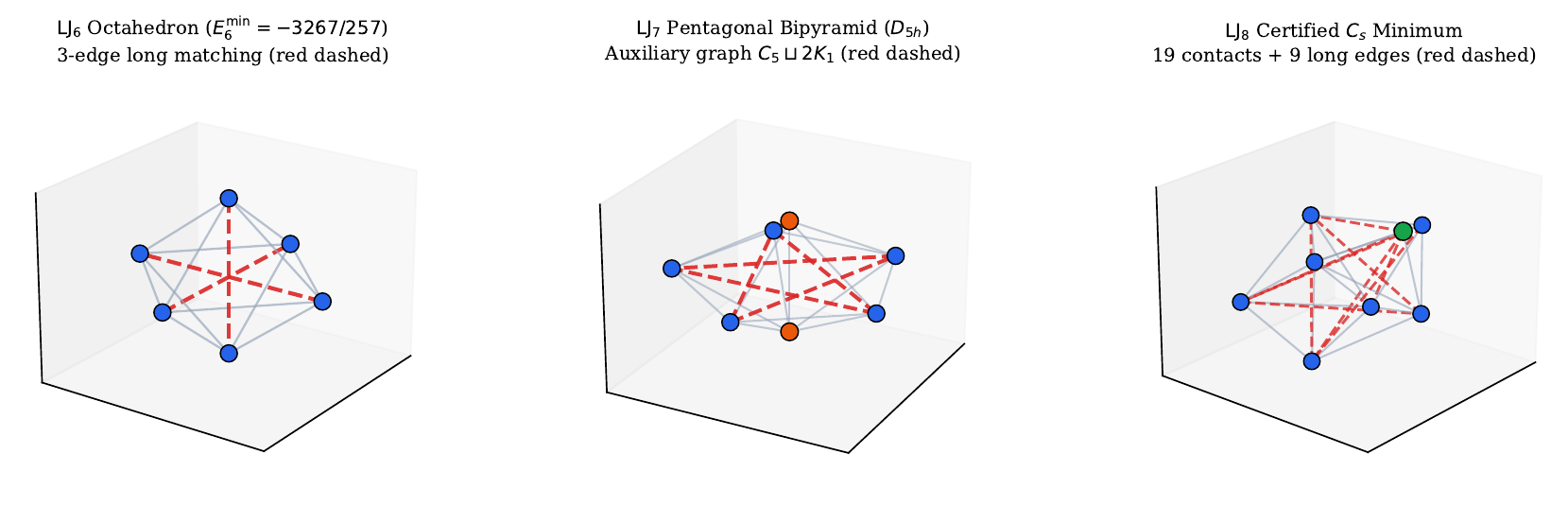}
\caption{Certified ground-state geometries and their surviving long/auxiliary distance graphs (gray solid lines: near-equilibrium contact edges; red dashed lines: long or auxiliary edges). Left: $\mathrm{LJ}_6$ regular octahedron with its $3$-edge antipodal matching. Center: $\mathrm{LJ}_7$ regular pentagonal bipyramid ($D_{5h}$), where the five equatorial diagonals form the auxiliary $5$-cycle $C_5$ and the two axial apices (orange) form $2K_1$. Right: Certified $\mathrm{LJ}_8$ $C_s$ local minimum (capped pentagonal bipyramid, capping vertex in green) with $19$ short contacts and $9$ long edges.}
\label{fig:clusters-graphs}
\end{figure}

\section{Exact Spectral Gram Relaxations and Global-to-Local Gluing}
\label{sec:spectral-framework}

Once the combinatorial topology of long and auxiliary edges is fixed, we must either prove that $E_N(X) > U$ on a squared-distance cell $\mathcal{B} = \prod_{(i,j)} [\ell_{ij}, u_{ij}]$ or localize every surviving configuration in $\mathcal{B}$ into a small neighborhood of the candidate minimizer. This section develops our general spectral-Gram relaxation and global-to-local gluing theory.

\subsection{Rational Bernstein Quadratic Minorants and Secant Cuts}
Let $f_\lambda(q) = (1+\lambda)q^{-6} - (2+\lambda)q^{-3}$ with $\lambda \in \Q$ ($\lambda = 0$ when stationarity is not invoked). Given a common leading curvature $K \in \Q_{>0}$, we construct for each pair $(i,j)$ on its interval $I_{ij} = [\ell_{ij}, u_{ij}] \subset (0,\infty)$ a quadratic lower support
\begin{equation}
\label{eq:scalar-quadratic-support}
f_\lambda(q) \ge K q^2 + b_{ij} q + c_{ij} \qquad \forall\, q \in I_{ij},
\end{equation}
with rational coefficients $b_{ij}, c_{ij} \in \Q$, certified by two exact algebraic mechanisms:
\begin{enumerate}[leftmargin=2em]
    \item \textbf{Short and intermediate intervals via exact Bernstein positivity:} Multiplying $f_\lambda(q) - (Kq^2 + bq + c) \ge 0$ by $q^6 > 0$ yields a degree-$8$ polynomial $P_8(q) = -Kq^8 - bq^7 - cq^6 - (2+\lambda)q^3 + (1+\lambda)$. When $Kq^2 + bq + c$ is constructed to be tangent to $f_\lambda$ at a rational contact point $z \in I_{ij}$, $P_8(q)$ factors exactly in $\Q[q]$ as $(q-z)^2 R_6(q)$. Transforming $R_6(q)$ to the Bernstein basis on $[\ell_{ij}, u_{ij}]$ (or a finite rational subdivision thereof), nonnegativity of all rational Bernstein coefficients certifies~\eqref{eq:scalar-quadratic-support} without floating-point approximation.
    \item \textbf{Long-distance intervals via concavity secants and interval multipliers:} For $q \ge B_0$ where $g_\lambda(q) := f_\lambda(q) - Kq^2$ is strictly concave ($g_\lambda''(q) < 0$), $g_\lambda(q)$ is bounded below on $[\ell_{ij}, u_{ij}]$ by the secant line connecting $(\ell_{ij}, g_\lambda(\ell_{ij}))$ and $(u_{ij}, g_\lambda(u_{ij}))$. In addition, for any nonnegative dual multipliers $\alpha_{ij}, \beta_{ij} \in \Q_{\ge 0}$, adding the nonpositive linear term $\alpha_{ij}(q_{ij} - u_{ij}) + \beta_{ij}(\ell_{ij} - q_{ij}) \le 0$ preserves validity of the lower bound on $[\ell_{ij}, u_{ij}]$.
\end{enumerate}

Summing~\eqref{eq:scalar-quadratic-support} over all $1 \le i < j \le N$ and writing $M_{ij} = M_{ji} := -b_{ij}/(2K)$ for $i \ne j$ and $M_{ii} := 0$ yields the general virial-shifted inequality $E_N(X) + \lambda V(X) \ge C_0 + K \sum_{i<j}(q_{ij} - M_{ij})^2$ for every $X \in \mathcal{B}$. Consequently, by Proposition~\ref{prop:scale-identity}, for every configuration $X \in \mathcal{B}$ when $\lambda = 0$, and for every configuration $X \in \mathcal{B}$ satisfying the virial condition $V(X) = 0$ (in particular, every stationary configuration $\nabla E_N(X) = 0$) when $\lambda \ne 0$, we obtain the completed-square distance lower bound
\begin{equation}
\label{eq:cell-quadratic-lower}
E_N(X) \ge C_0 + K \sum_{1 \le i < j \le N} (q_{ij} - M_{ij})^2, \qquad C_0 := \sum_{1 \le i < j \le N} \left(c_{ij} - \frac{b_{ij}^2}{4K}\right).
\end{equation}

\subsection{Centered Gram Matrix Identity and Dual Square Completion}
For any configuration $X = (x_1,\dots,x_N) \in \Omega_N$ centered at the origin ($\sum_{i=1}^N x_i = 0$), let $G \in \R^{N \times N}$ denote its Gram matrix $G_{ij} = \langle x_i, x_j \rangle$. Then
\begin{equation}
\label{eq:gram-constraints}
G \succeq 0, \qquad G\One = 0, \qquad \rank(G) \le 3,
\end{equation}
and the squared distances are recovered via $q_{ij} = G_{ii} + G_{jj} - 2G_{ij}$. Conversely, if $Q \in \R^{N \times N}$ is the hollow symmetric squared-distance matrix ($Q_{ij} = q_{ij}$, $Q_{ii} = 0$) and $J := I_N - \frac{1}{N}\One\One^\top$ is the centering projector, then $G = -\frac{1}{2} J Q J$.

\begin{proposition}[Exact Frobenius Identity and Dual Square Completion]
\label{prop:gram-identity}
Let $M \in \R^{N \times N}$ be any symmetric hollow matrix ($M_{ii} = 0$) and define $G_0 := -\frac{1}{2} J M J$. Then for every centered Gram matrix $G$ with induced squared distances $q_{ij} = G_{ii} + G_{jj} - 2G_{ij}$, the matrix difference $A := G - G_0$ satisfies the exact identity
\begin{equation}
\label{eq:frobenius-distance-identity}
\sum_{1 \le i < j \le N} (q_{ij} - M_{ij})^2 = 2\|A\|_F^2 + N\|\diag(A)\|_2^2 + \bigl(\tr(A)\bigr)^2.
\end{equation}
Furthermore, for $N \ge 3$, any dual vector $v \in \Q^N$, and any scalar $\tau \in \Q$, define the symmetric matrices
\begin{equation}
\label{eq:W-and-H}
W := J\bigl(N\diag(v) + \tau I_N\bigr)J, \qquad H := G_0 - \frac{1}{2}W.
\end{equation}
Then every configuration $X \in \mathcal{B}$ satisfying~\eqref{eq:cell-quadratic-lower} (namely, every $X \in \mathcal{B}$ when $\lambda = 0$, or every $X \in \mathcal{B}$ with $V(X) = 0$ when $\lambda \ne 0$) obeys
\begin{equation}
\label{eq:completed-square-gram}
E_N(X) \ge \mathcal{B}_{\mathrm{spec}}(H) + 2K\bigl(\|G - H\|_F^2 - \|T - H\|_F^2\bigr) + K\Bigl(N\|\diag(A) - v\|_2^2 + (\tr(A) - \tau)^2\Bigr),
\end{equation}
where $T$ is the positive top-three spectral truncation of $H$, and
\begin{equation}
\label{eq:B-spec}
\mathcal{B}_{\mathrm{spec}}(H) := C_0 + K\left(2\|H\|_F^2 - 2\sum_{r=1}^3 \bigl(\max(\lambda_r(H), 0)\bigr)^2 - \frac{1}{2}\|W\|_F^2 - N\|v\|_2^2 - \tau^2\right),
\end{equation}
with $\lambda_1(H) \ge \lambda_2(H) \ge \dots \ge \lambda_N(H)$ denoting the eigenvalues of $H$ in nonincreasing order.
\end{proposition}

\begin{proof}
Since $J\One = 0$, $G_0\One = 0$ and $A\One = 0$. From $q_{ij} = G_{ii} + G_{jj} - 2G_{ij}$ and $M_{ij} = (G_0)_{ii} + (G_0)_{jj} - 2(G_0)_{ij}$ (which holds because $M$ is hollow and $G_0 = -\frac{1}{2}JMJ$), the residuals $d_{ij} := q_{ij} - M_{ij}$ satisfy $d_{ij} = A_{ii} + A_{jj} - 2A_{ij}$ for $i \ne j$ and $d_{ii} = 0$. Expanding $\frac{1}{2}\sum_{i,j=1}^N (A_{ii} + A_{jj} - 2A_{ij})^2$ using $\sum_j A_{ij} = 0$ yields
\[
\sum_{1 \le i < j \le N} d_{ij}^2 = \frac{1}{2}\sum_{i,j=1}^N \bigl( (A_{ii}+A_{jj})^2 - 4(A_{ii}+A_{jj})A_{ij} + 4A_{ij}^2 \bigr) = N\sum_{i=1}^N A_{ii}^2 + \left(\sum_{i=1}^N A_{ii}\right)^2 + 2\sum_{i,j=1}^N A_{ij}^2,
\]
which proves~\eqref{eq:frobenius-distance-identity}. Next, since $A = JAJ$, we have
\[
N\|\diag(A)\|_2^2 + (\tr(A))^2 = N\|\diag(A)-v\|_2^2 + (\tr(A)-\tau)^2 + 2\langle A, W \rangle - N\|v\|_2^2 - \tau^2.
\]
Combining this with $2\|A\|_F^2 + 2\langle A, W \rangle = 2\|G - H\|_F^2 - \frac{1}{2}\|W\|_F^2$ and noting that $\|T - H\|_F^2 = \|H\|_F^2 - \sum_{r=1}^3 (\max(\lambda_r(H),0))^2$ establishes~\eqref{eq:completed-square-gram} and~\eqref{eq:B-spec}.
\end{proof}

\begin{remark}[Exact Certification of Eigenvalue Bounds via $LDL^\top$ Inertia]
\label{rem:ldl-inertia}
By the Eckart--Young theorem~\cite{eckart1936approximation} for positive semidefinite matrices of rank at most three, $\|G - H\|_F^2 \ge \|T - H\|_F^2$ holds for all $G \succeq 0$ with $\rank(G) \le 3$. Thus $\mathcal{B}_{\mathrm{spec}}(H)$ in~\eqref{eq:B-spec} is a valid energy lower bound under~\eqref{eq:cell-quadratic-lower}. To certify $\mathcal{B}_{\mathrm{spec}}(H)$ in exact arithmetic without floating-point eigensolvers, we take rational upper bounds $s_1 \ge s_2 \ge s_3 \ge 0$ and verify $\lambda_r(H) \le s_r$ by computing the exact rational $LDL^\top$ factorization of $s_r I_N - H$ and counting its negative pivots via Sylvester's law of inertia.
\end{remark}

\subsection{The Rank-\texorpdfstring{$d$}{d} Frobenius Spectral-Gap Localization Lemma}
When a distance cell $\mathcal{B}$ contains or abuts the global minimizer, $\mathcal{B}_{\mathrm{spec}}(H)$ cannot exceed $U$. Instead, for any ambient dimension $1 \le d < N$ (with $d=3$ and $N \ge 4$ in our three-dimensional applications), we exploit the spectral gap between $\lambda_d(H)$ and $\lambda_{d+1}(H)$ to bound the Frobenius distance $\|G - T\|_F$.

\begin{lemma}[Quadratic Growth of Rank-$d$ Spectral Approximation]
\label{lem:spectral-gap}
Fix integers $1 \le d < N$, and let $H \in \R^{N \times N}$ be symmetric with eigenvalues $\lambda_1(H) \ge \lambda_2(H) \ge \dots \ge \lambda_N(H)$. Suppose there exist real numbers $a, b$ (chosen in $\Q$ for exact certification) such that
\begin{equation}
\label{eq:eigengap-assumption}
\lambda_d(H) \ge a > b \ge \max\bigl(\lambda_{d+1}(H),\, 0\bigr), \qquad c := 1 - \frac{b}{a} \in (0,1].
\end{equation}
Let $T = \sum_{r=1}^d \lambda_r(H) u_r u_r^\top$ be the positive top-$d$ spectral truncation of $H$. Then for every positive semidefinite matrix $G \in \R^{N \times N}$ with $\rank(G) \le d$,
\begin{equation}
\label{eq:spectral-gap-inequality}
\|G - H\|_F^2 - \|T - H\|_F^2 \ge c\,\|G - T\|_F^2,
\end{equation}
and the constant $c = 1 - b/a$ is uniformly sharp over this class of matrices $(H,G)$ for any given $a > b \ge 0$.
\end{lemma}

\begin{proof}
Let $P_{\mathrm{top}} = \sum_{r=1}^d u_r u_r^\top$ be the orthogonal projector onto the top $d$ eigenspaces of $H$, and let $P_\perp = I_N - P_{\mathrm{top}}$. Write $H_\perp := H - T$. Since $H_\perp P_{\mathrm{top}} = 0$ and $\lambda_{d+1}(H) \le b$ with $b \ge 0$, we have $H_\perp \preceq bP_\perp$. Because $G \succeq 0$, it follows that $\langle H_\perp, G \rangle \le b\tr(P_\perp G)$. Expanding $\|G - H\|_F^2 = \|(G - T) - H_\perp\|_F^2$ and using $\langle T, H_\perp \rangle = 0$ gives
\begin{equation}
\label{eq:gap-step1}
\|G - H\|_F^2 - \|T - H\|_F^2 = \|G - T\|_F^2 - 2\langle H_\perp, G \rangle \ge \|G - T\|_F^2 - 2b\tr(P_\perp G).
\end{equation}
Next, consider the symmetric matrix $\widetilde{H} := T + aP_\perp$, whose $d$ largest eigenvalues are $\lambda_1(H) \ge \dots \ge \lambda_d(H) \ge a$ (with spectral truncation $T$) and whose remaining $N-d$ eigenvalues equal $a$ (on $\operatorname{im} P_\perp$). By the Eckart--Young best rank-$d$ approximation theorem~\cite{eckart1936approximation}, $T$ minimizes $\|M - \widetilde{H}\|_F^2$ over all matrices $M$ of rank at most $d$. Since $\rank(G) \le d$, we have $\|G - (T + aP_\perp)\|_F^2 \ge \|T - (T + aP_\perp)\|_F^2 = \|aP_\perp\|_F^2$. Expanding the left-hand side using $\langle T, P_\perp \rangle = 0$ yields
\begin{equation}
\label{eq:gap-step2}
\|G - T\|_F^2 - 2a\tr(P_\perp G) + \|aP_\perp\|_F^2 \ge \|aP_\perp\|_F^2 \implies \|G - T\|_F^2 \ge 2a\tr(P_\perp G).
\end{equation}
Since $b \ge 0$ and $a > 0$, combining~\eqref{eq:gap-step1} and~\eqref{eq:gap-step2} gives $\|G - H\|_F^2 - \|T - H\|_F^2 \ge (1 - b/a)\|G - T\|_F^2$. For sharpness, take
\[
H = \diag(aI_d,b,0_{N-d-1}), \qquad
T = \diag(aI_d,0_{N-d}), \qquad
G = \diag(aI_{d-1},0,a,0_{N-d-1}),
\]
where $0_m$ denotes the $m \times m$ zero matrix and blocks of size zero are omitted. Then $\rank(G)=d$ and $\|G - H\|_F^2 - \|T - H\|_F^2 = 2a(a-b) = (1-b/a)\|G - T\|_F^2$.
\end{proof}

Combining Proposition~\ref{prop:gram-identity} and Lemma~\ref{lem:spectral-gap} (with $d=3$) bounds the Frobenius distance between any low-energy Gram matrix $G$ and the exact candidate Gram matrix $G_{\mathrm{cand}}$ without computing the irrational eigenvectors of $T$.

\begin{corollary}[Certificate-Checkable Localization to the Candidate Gram Matrix]
\label{cor:gram-localization}
Let $N \ge 4$ and $d = 3$. Under the hypotheses of Proposition~\ref{prop:gram-identity} and Lemma~\ref{lem:spectral-gap}, let $\mathrm{LB} \le \mathcal{B}_{\mathrm{spec}}(H)$ be a certified rational lower bound on a cell $\mathcal{B}$, and let $G_{\mathrm{cand}} \succeq 0$ be the exact rank-three centered Gram matrix of the candidate minimizer. Suppose $s_1,s_2,s_3 \in \Q_{>0}$ satisfy $\lambda_r(H) \le s_r$ ($r=1,2,3$), and let $\delta \in \Q_{\ge 0}$ be a certified upper bound on
\begin{equation}
\label{eq:delta-bound}
\|G_{\mathrm{cand}}\|_F^2 - 2\langle G_{\mathrm{cand}}, H \rangle + \sum_{r=1}^3 s_r^2 \le \delta.
\end{equation}
Consider any configuration $X \in \mathcal{B}$ with $E_N(X) \le U$ that satisfies~\eqref{eq:cell-quadratic-lower} (i.e., any $X \in \mathcal{B}$ with $E_N(X) \le U$ when $\lambda = 0$, or any $X \in \mathcal{B}$ with $E_N(X) \le U$ and $V(X) = 0$ when $\lambda \ne 0$). If $\mathrm{LB} > U$, no such configuration exists and $\mathcal{B}$ is excluded. If $\mathrm{LB} \le U$, then for any chosen rational upper bound $\eta \in \Q_{>0}$ satisfying $\sqrt{(U - \mathrm{LB})/(2Kc)} + \sqrt{\delta/c} \le \eta$, the centered Gram matrix $G$ of $X$ obeys
\begin{equation}
\label{eq:gram-cand-bound}
\|G - G_{\mathrm{cand}}\|_F \le \sqrt{\frac{U - \mathrm{LB}}{2Kc}} + \sqrt{\frac{\delta}{c}} \le \eta,
\end{equation}
and every pairwise squared distance $q_{ij}$ of $X$ is enclosed within
\begin{equation}
\label{eq:pair-distance-localization}
|q_{ij} - q_{ij}^{\mathrm{cand}}| \le 2\eta \qquad \forall\, 1 \le i < j \le N.
\end{equation}
\end{corollary}

\begin{proof}
By~\eqref{eq:completed-square-gram} and Lemma~\ref{lem:spectral-gap}, $E_N(X) \ge \mathrm{LB} + 2Kc\|G - T\|_F^2$. Thus $E_N(X) \le U$ is impossible when $\mathrm{LB} > U$, and implies $\|G - T\|_F \le \sqrt{(U - \mathrm{LB})/(2Kc)}$ when $\mathrm{LB} \le U$. Similarly, since $\rank(G_{\mathrm{cand}}) \le 3$ and $\|T - H\|_F^2 = \|H\|_F^2 - \sum_{r=1}^3 \lambda_r(H)^2 \ge \|H\|_F^2 - \sum_{r=1}^3 s_r^2$, Lemma~\ref{lem:spectral-gap} applied to $G_{\mathrm{cand}}$ gives
\[
c\|G_{\mathrm{cand}} - T\|_F^2 \le \|G_{\mathrm{cand}} - H\|_F^2 - \|T - H\|_F^2 \le \|G_{\mathrm{cand}}\|_F^2 - 2\langle G_{\mathrm{cand}}, H \rangle + \sum_{r=1}^3 s_r^2 \le \delta.
\]
The triangle inequality $\|G - G_{\mathrm{cand}}\|_F \le \|G - T\|_F + \|T - G_{\mathrm{cand}}\|_F$ then yields~\eqref{eq:gram-cand-bound}. Finally, for $w_{ij} := e_i - e_j$, we have $q_{ij} - q_{ij}^{\mathrm{cand}} = \langle w_{ij}w_{ij}^\top, G - G_{\mathrm{cand}} \rangle$ and $\|w_{ij}w_{ij}^\top\|_F = 2$, so Cauchy--Schwarz gives~\eqref{eq:pair-distance-localization}.
\end{proof}

\subsection{Oriented Cartesian Reconstruction and Interval Newton Gluing}
\label{subsec:newton-gluing}
Given the narrow distance intervals $q_{ij} \in [q_{ij}^{\mathrm{cand}} - 2\eta,\, q_{ij}^{\mathrm{cand}} + 2\eta]$ from Corollary~\ref{cor:gram-localization}, we reconstruct a full $(3N-6)$-dimensional gauge-fixed Cartesian coordinate box and contract it into the local uniqueness neighborhood of the candidate minimizer.

To make all matrix dimensions and domain requirements explicit, let $d_{\mathrm{red}} := 3N-6$ and let $\iota: \R^{d_{\mathrm{red}}} \to (\R^3)^N$ denote the linear gauge-fixing embedding that places particle $0$ at $p_0 = (0,0,0)$, particle $1$ at $p_1 = (d,0,0)$, particle $2$ at $p_2 = (x_2,y_2,0)$, and particles $i = 3,\dots,N-1$ at $p_i = (x_i,y_i,z_i)$ for the reduced coordinate vector $x = (d, x_2, y_2, x_3, y_3, z_3, \dots, x_{N-1}, y_{N-1}, z_{N-1}) \in \R^{d_{\mathrm{red}}}$. On the open collision-free domain $\mathcal{U} := \iota^{-1}(\Omega_N) \subset \R^{d_{\mathrm{red}}}$, define the smooth reduced energy $F := E_N \circ \iota : \mathcal{U} \to \R$, with reduced gradient $\nabla F(x) \in \R^{d_{\mathrm{red}}}$ and reduced Hessian $\nabla^2 F(x) \in \R^{d_{\mathrm{red}} \times d_{\mathrm{red}}}$.

\begin{remark}[Reduced Gauge-Fixed Stationarity and Hessian]
\label{rem:reduced-hessian}
On the non-collinear open domain $\{x \in \mathcal{U} : d > 0,\, y_2 > 0\}$, the gauge subspace $\operatorname{im}(d\iota) \subset (\R^3)^N$ and the six-dimensional $\mathrm{SE}(3)$-orbit tangent space $T_{\iota(x)}\bigl(\mathrm{SE}(3)\cdot\iota(x)\bigr)$ are complementary. Indeed, any infinitesimal rigid motion $\xi_i = a + \omega \times p_i$ ($a,\omega \in \R^3$) evaluated on the six fixed gauge coordinates $(p_{0x},p_{0y},p_{0z},p_{1y},p_{1z},p_{2z})$ yields $(a_x, a_y, a_z, a_y + \omega_z d, a_z - \omega_y d, a_z + \omega_x y_2 - \omega_y x_2)$, a lower-triangular linear system in $(a_x, a_y, a_z, \omega_z, \omega_y, \omega_x)$ of determinant $-d^2 y_2 \ne 0$. Thus $\operatorname{im}(d\iota) \cap T_{\iota(x)}\bigl(\mathrm{SE}(3)\cdot\iota(x)\bigr) = \{0\}$ and $(\R^3)^N = \operatorname{im}(d\iota) \oplus T_{\iota(x)}\bigl(\mathrm{SE}(3)\cdot\iota(x)\bigr)$. Because $E_N$ is $\mathrm{SE}(3)$-invariant, $\langle \nabla E_N(\iota(x)), \xi \rangle = 0$ for all $\xi \in T_{\iota(x)}\bigl(\mathrm{SE}(3)\cdot\iota(x)\bigr)$, so reduced stationarity $\nabla F(x^*) = d\iota^\top \nabla E_N(\iota(x^*)) = 0$ is equivalent to full Cartesian stationarity $\nabla E_N(\iota(x^*)) = 0$. Furthermore, positive definiteness of $\nabla^2 F \succ 0$ on a convex neighborhood of $x^*$ in $\{d > 0,\, y_2 > 0\}$ implies nondegenerate strict local minimality of $E_N$ modulo $\mathrm{E}(3)$. Note that while $d > 0$ and $y_2 > 0$ fix the continuous $\mathrm{SE}(3)$ gauge, the half-space condition $z_3 \ge 0$ selects a representative modulo the mirror reflection $z \mapsto -z$ across the $xy$-plane only when $z_3 > 0$ (if $z_3 = 0$, $p_0,p_1,p_2,p_3$ are coplanar and the reflection $(z_4,\dots,z_{N-1}) \mapsto -(z_4,\dots,z_{N-1})$ is not fixed by $z_3 \ge 0$). In Proposition~\ref{prop:chart-newton}, the strict non-coplanarity condition $z_3 > 0$ both uniquely fixes this reflection and supplies the positive denominator required to reconstruct $z_i$ ($i \ge 4$).
\end{remark}

\begin{proposition}[Oriented Gauge Chart Reconstruction and Interval Newton Gluing]
\label{prop:chart-newton}
Suppose that on a distance cell $\mathcal{B}$ under threshold $U$, every stationary configuration $X \in \mathcal{B}$ with $E_N(X) \le U$ has its squared distances (up to an admissible particle relabeling in $S_N$, such as the uniform enclosure $[q_{ij}] = [q_{ij}^{\mathrm{cand}} - 2\eta, q_{ij}^{\mathrm{cand}} + 2\eta]$ supplied by Corollary~\ref{cor:gram-localization}) enclosed in outward-rounded intervals $[q_{ij}]$. Here and below, for any scalar, vector, or matrix expression $\phi$, we denote by $[\phi]$ its outward-rounded interval extension over the enclosing intervals of its arguments; a symmetric interval matrix $\mathcal{A} = ([A_{ij}^{\mathrm{lo}}, A_{ij}^{\mathrm{up}}])$ is said to be positive definite ($\mathcal{A} \succ 0$) if every real symmetric matrix $A \in \mathcal{A}$ is positive definite. By an isometry in $\mathrm{E}(3)$, represent any such $X$ in the gauge chart $\iota(x)$ with $d > 0$, $y_2 > 0$, and $z_3 \ge 0$, and set $s_{ij} := \frac{1}{2}(q_{0i} + q_{0j} - q_{ij})$. If the interval enclosures defined sequentially from $[q_{ij}]$ and $[s_{ij}] := \frac{1}{2}([q_{0i}] + [q_{0j}] - [q_{ij}])$ verify strict positivity $\inf[q_{01}] > 0$, $\inf([q_{02}] - [x_2]^2) > 0$, and $\inf([q_{03}] - [x_3]^2 - [y_3]^2) > 0$ (so that $d > 0$, $y_2 > 0$, and $z_3 > 0$, uniquely fixing the reflection across the $xy$-plane), then for every stationary configuration $X \in \mathcal{B}$ with $E_N(X) \le U$, its gauge coordinate vector $x \in \R^{d_{\mathrm{red}}}$ ($d_{\mathrm{red}} = 3N-6$) is enclosed in the single compact coordinate box $\mathcal{X}_0 := [d] \times [x_2] \times [y_2] \times \prod_{i=3}^{N-1}([x_i]\times [y_i]\times [z_i]) \subset \R^{d_{\mathrm{red}}}$ reconstructed via
\begin{equation}
\label{eq:chart-reconstruction}
\begin{aligned}
d &= \sqrt{q_{01}}, \qquad & x_i &= \frac{s_{1i}}{d} \quad (i \ge 2), \\
y_2 &= \sqrt{q_{02} - x_2^2}, \qquad & y_i &= \frac{s_{2i} - x_2 x_i}{y_2} \quad (i \ge 3), \\
z_3 &= \sqrt{q_{03} - x_3^2 - y_3^2}, \qquad & z_i &= \frac{s_{3i} - x_3 x_i - y_3 y_i}{z_3} \quad (i \ge 4).
\end{aligned}
\end{equation}
For any collision-free box $\mathcal{X} \subset \mathcal{U}$, chosen center $\hat{x} \in \mathcal{X}$, rational preconditioner $Y \in \Q^{d_{\mathrm{red}} \times d_{\mathrm{red}}}$, and any interval matrix enclosure $\mathcal{M}_Y(\mathcal{X}) \supseteq \{Y \nabla^2 F(\xi) : \xi \in \mathcal{X}\}$ of the preconditioned reduced Hessian on $\mathcal{X}$ (obtained either as the interval product $Y[\nabla^2 F(\mathcal{X})]$ or as a direct pair-correlated preconditioned enclosure $[Y \nabla^2 F(\mathcal{X})]$), define the raw Krawczyk image and the root-retaining Krawczyk--Newton operator on $F = E_N \circ \iota$ by
\begin{equation}
\label{eq:krawczyk-step}
\mathcal{K}(\mathcal{X}) := \hat{x} - Y \nabla F(\hat{x}) + \bigl(I_{d_{\mathrm{red}}} - \mathcal{M}_Y(\mathcal{X})\bigr)(\mathcal{X} - \hat{x}), \qquad
\mathcal{N}(\mathcal{X}) := \mathcal{K}(\mathcal{X}) \cap \mathcal{X}.
\end{equation}
Suppose further that:
\begin{enumerate}[label=(\roman*),leftmargin=2.2em]
    \item \textbf{Local existence, uniqueness, and strict local minimality on $\mathcal{X}_{\mathrm{loc}}$:} There exist a collision-free target box $\mathcal{X}_{\mathrm{loc}} \subset \mathcal{U}$ with center $\hat{x}_{\mathrm{loc}} \in \mathcal{X}_{\mathrm{loc}}$ and an invertible upper-triangular rational preconditioner $P \in \mathrm{GL}_{d_{\mathrm{red}}}(\Q)$ with $Y_{\mathrm{loc}} = P P^\top \in \mathrm{GL}_{d_{\mathrm{red}}}(\Q)$ such that the \emph{unintersected} Krawczyk image satisfies $\mathcal{K}_{\mathrm{loc}}(\mathcal{X}_{\mathrm{loc}}) \subset \operatorname{int}(\mathcal{X}_{\mathrm{loc}})$, $\|I_{d_{\mathrm{red}}} - \mathcal{M}_{Y_{\mathrm{loc}}}(\mathcal{X}_{\mathrm{loc}})\|_\infty < 1$, and the preconditioned symmetric interval Hessian $\mathcal{H}_{\mathrm{pre}} := P^\top [\nabla^2 F(\mathcal{X}_{\mathrm{loc}})] P \succ 0$ (certified via Gershgorin's circle theorem by verifying $\min_{1 \le i \le d_{\mathrm{red}}}\bigl((\mathcal{H}_{\mathrm{pre}})_{ii}^{\mathrm{lo}} - \sum_{j \ne i}|(\mathcal{H}_{\mathrm{pre}})_{ij}|^{\mathrm{up}}\bigr) > 0$; or, when an exact stationary point $x^* \in \operatorname{int}(\mathcal{X}_{\mathrm{loc}})$ is already algebraically constructed, $\mathcal{H}_{\mathrm{pre}} \succ 0$ on the convex box $\mathcal{X}_{\mathrm{loc}}$).
    \item \textbf{Full-box non-collision on $\mathcal{X}_0$ and global-to-local root retention into $\mathcal{X}_{\mathrm{loc}}$:} Outward interval evaluation on the Cartesian box $\mathcal{X}_0$ verifies $\inf_{x \in \mathcal{X}_0}\|p_i(x) - p_j(x)\|_2^2 > 0$ for all $0 \le i < j \le N-1$ (so $\mathcal{X}_0 \subset \mathcal{U}$ is collision-free), and starting from $\mathcal{X}_0$, $m \ge 1$ iterations $\mathcal{X}_{k} = \mathcal{N}_k(\mathcal{X}_{k-1})$ (each evaluated with a center $\hat{x}_{k-1} \in \mathcal{X}_{k-1}$ and a full collision-free box enclosure $\mathcal{M}_{Y_{k-1}}(\mathcal{X}_{k-1})$ on $\mathcal{X}_{k-1} \subseteq \mathcal{X}_0 \subset \mathcal{U}$) satisfy $\mathcal{X}_m \subset \mathcal{X}_{\mathrm{loc}}$.
\end{enumerate}
Then $\mathcal{X}_{\mathrm{loc}}$ contains a unique stationary point $x^*$ of $F$, which is a strict local minimizer of $E_N$ modulo rigid motions $\mathrm{E}(3)$, and every stationary configuration in the distance cell $\mathcal{B}$ with $E_N(X) \le U$ coincides (up to $\mathrm{E}(3)$ and $S_N$) with $\iota(x^*)$.
\end{proposition}

\begin{proof}
For any stationary configuration $X \in \mathcal{B}$ with $E_N(X) \le U$, its relabeled squared distances lie in $[q_{ij}]$ by hypothesis. Since $s_{ij} = \langle p_i, p_j \rangle$, the strict lower bounds $\inf[q_{01}] > 0$, $\inf([q_{02}]-[x_2]^2) > 0$, and $\inf([q_{03}]-[x_3]^2-[y_3]^2) > 0$ force $d, y_2, z_3 > 0$ under $z_3 \ge 0$, uniquely fixing the reflection across the $xy$-plane. Evaluating inner products successively with $p_1 = (d,0,0)$, $p_2 = (x_2,y_2,0)$, and $p_3 = (x_3,y_3,z_3)$ yields~\eqref{eq:chart-reconstruction}; for $i \ge 4$, $z_i$ is recovered linearly from $s_{3i}$ using the single positive denominator $z_3 > 0$, preserving all relative $z$-signs and placing $x \in \mathcal{X}_0$. Because independent Cartesian coordinate intervals can lose pairwise distance correlations, verifying $\inf_{x \in \mathcal{X}_0}\|p_i(x) - p_j(x)\|_2^2 > 0$ for all $i < j$ ensures $\mathcal{X}_0 \subset \mathcal{U}$ (and hence $\mathcal{X}_k \subseteq \mathcal{X}_0 \subset \mathcal{U}$ for all $k \ge 1$) so $F$ is smooth on every box $\mathcal{X}_k$. Next, on $\mathcal{X}_{\mathrm{loc}}$, the unintersected inclusion $\mathcal{K}_{\mathrm{loc}}(\mathcal{X}_{\mathrm{loc}}) \subset \operatorname{int}(\mathcal{X}_{\mathrm{loc}})$ together with $\|I_{d_{\mathrm{red}}} - \mathcal{M}_{Y_{\mathrm{loc}}}(\mathcal{X}_{\mathrm{loc}})\|_\infty < 1$ and invertibility of $Y_{\mathrm{loc}}$ implies by the Krawczyk--Banach contraction theorem that $F$ has a unique stationary point $x^* \in \operatorname{int}(\mathcal{X}_{\mathrm{loc}})$, which satisfies $\nabla E_N(\iota(x^*)) = 0$ by Remark~\ref{rem:reduced-hessian}. Moreover, for every real symmetric matrix $H \in [\nabla^2 F(\mathcal{X}_{\mathrm{loc}})]$, the symmetric matrix $M = P^\top H P \in \mathcal{H}_{\mathrm{pre}}$ has minimum eigenvalue bounded below by $\min_i\bigl(M_{ii} - \sum_{j \ne i}|M_{ij}|\bigr) \ge \min_i\bigl((\mathcal{H}_{\mathrm{pre}})_{ii}^{\mathrm{lo}} - \sum_{j \ne i}|(\mathcal{H}_{\mathrm{pre}})_{ij}|^{\mathrm{up}}\bigr) > 0$ by Gershgorin's theorem; since $P \in \mathrm{GL}_{d_{\mathrm{red}}}(\Q)$, Sylvester's law of inertia implies $H \succ 0$ and hence $\nabla^2 F(\xi) \succ 0$ throughout the convex box $\mathcal{X}_{\mathrm{loc}}$ (which also directly yields uniqueness and strict local minimality modulo $\mathrm{E}(3)$ when $x^* \in \operatorname{int}(\mathcal{X}_{\mathrm{loc}})$ is constructed algebraically). Finally, for any convex box $\mathcal{X} \subset \mathcal{U}$ with $\hat{x} \in \mathcal{X}$ and any stationary point $\tilde{x} \in \mathcal{X}$ ($\nabla F(\tilde{x}) = 0$), the segment $\{\hat{x} + \theta(\tilde{x}-\hat{x}) : \theta \in [0,1]\}$ lies in $\mathcal{X}$; since interval matrices are convex, $\int_0^1 Y\nabla^2 F(\hat{x} + \theta(\tilde{x}-\hat{x}))\,d\theta \in \mathcal{M}_Y(\mathcal{X})$, so the mean value theorem gives
\[
\tilde{x} = \hat{x} - Y\nabla F(\hat{x}) + \left(I_{d_{\mathrm{red}}} - \int_0^1 Y\nabla^2 F\bigl(\hat{x} + \theta(\tilde{x}-\hat{x})\bigr)\,d\theta\right)(\tilde{x} - \hat{x}) \in \mathcal{K}(\mathcal{X}) \cap \mathcal{X} = \mathcal{N}(\mathcal{X}).
\]
By induction over $k=1,\dots,m$, every stationary point in $\mathcal{X}_0$ lies in $\mathcal{X}_m \subset \mathcal{X}_{\mathrm{loc}}$ and therefore coincides with $x^*$.
\end{proof}

\section{Complete Proof of Global Optimality for \texorpdfstring{$N=6$}{N=6}}
\label{sec:lj6}

In this section we prove Theorem~\ref{thm:lj6}, establishing that the global minimum of $E_6$ on $\Omega_6$ is the exact rational number $-3267/257$, attained uniquely by the regular octahedron.

\subsection{The Octahedral Scale Identity and Reduction to Matchings}
Place six particles at the vertices of a regular octahedron $\pm a e_1, \pm a e_2, \pm a e_3 \in \R^3$ ($a > 0$), and let $q = 2a^2$ be the squared length of its twelve short edges. The three antipodal pairs have squared distance $4a^2 = 2q$. Evaluating~\eqref{eq:lj-energy} gives the exact rational completion of squares:
\begin{equation}
\label{eq:octahedron-energy}
E_6^{\mathrm{oct}}(q) = 12(q^{-6} - 2q^{-3}) + 3\bigl((2q)^{-6} - 2(2q)^{-3}\bigr) = \frac{771}{64}q^{-6} - \frac{99}{4}q^{-3} = \frac{771}{64}\left(q^{-3} - \frac{264}{257}\right)^2 - \frac{3267}{257}.
\end{equation}
Thus $E_6^{\mathrm{oct}}(q)$ attains its unique minimum $U_6 := -3267/257 \approx -12.7120622568$ at $q_*^3 = 257/264$. By octahedral symmetry, the full Cartesian gradient at $a e_1$ is $\nabla_{ae_1} E_6 = \frac{3a}{16q^7}(264q^3 - 257)e_1$, which vanishes at $q = q_*$. By Lemma~\ref{lem:attainment} and Remark~\ref{rem:elementary-attainment}, global minimizers of $E_6$ exist in $\mathcal{S}_6(U_6)$ and satisfy $\nabla E_6(X) = 0$.

Set the true-long distance threshold to $r_{\mathrm{long}} = 27/20 = 1.35$ ($B_0 = r_{\mathrm{long}}^2 = 729/400 = 1.8225$) and the auxiliary threshold to $C = 3/2 = 1.5$. Using the certified five-particle lower bound $L_5 = -71123847/7812500$ (Proposition~\ref{prop:lj5-foundation}), any vertex incident to two true-long edges would force
\[
E_6(X) \ge L_5 - 3 + 2v(27/20) = -12.71005426328\dots > -\frac{3267}{257} = U_6.
\]
By Proposition~\ref{prop:graph-pruning}(i), $G_{\mathrm{long}}(X)$ has maximum degree at most $1$, so it is a matching of size $k \in \{0,1,2,3\}$. Moreover, for any two particles $i,j$, each of $i$ and $j$ is incident to at most one true-long edge among the remaining four particles; thus at least two of the remaining four particles are connected to \emph{both} $i$ and $j$ by non-long edges ($< r_{\mathrm{long}}$), proving by Proposition~\ref{prop:graph-pruning}(iv) that $r_{ij} < 2r_{\mathrm{long}}$ ($q_{ij} < 4B_0 = 729/100$) for all $15$ pairs. Similarly, the inequalities $L_5 - 2 + 3f(3/2) > U_6$ and $-15 + 5\ell(3/2) > U_6$ imply that $G_{\mathrm{aux}}(X)$ has maximum degree at most $2$ and at most $4$ auxiliary edges, while $-15 + 3\ell(B_0) + \ell(3/2) \approx -12.41410387 > U_6$ rules out any fourth auxiliary edge when $k=3$, ensuring that for $k=3$, $G_{\mathrm{aux}}(X)$ coincides with the $3$-edge true-long matching.

\paragraph{Exclusion of $k=0$ (no true-long edges).}
When all $15$ squared distances satisfy $0 < q_{ij} \le B_0 = 729/400$, an exact rational Bernstein check establishes the global quadratic minorant ($\lambda = 0$)
\begin{equation}
f(q) \ge \frac{3}{2}q^2 - \frac{3381}{1000}q + \frac{219}{250} \qquad \forall\, q \in (0, 729/400].
\end{equation}
For any centered rank-at-most-three Gram matrix $G \in \R^{6 \times 6}$, writing $S_1 = \sum_{i<j} q_{ij} = 6\tr(G)$ and applying~\eqref{eq:frobenius-distance-identity} with $M=0$ yields $\sum_{i<j} q_{ij}^2 = 2\tr(G^2) + 6\sum_{i=1}^6 G_{ii}^2 + (\tr G)^2 \ge (2\cdot\frac{1}{3} + 6\cdot\frac{1}{6} + 1)(\tr G)^2 = \frac{8}{3}(\tr G)^2 = \frac{2}{27}S_1^2$. Minimizing the resulting scalar quadratic in $S_1$ excludes the entire $k=0$ domain $(0, B_0]^{15}$ in a single analytical step without enumerating auxiliary graphs:
\[
E_6(X) \ge \frac{3}{2}\cdot\frac{2}{27}S_1^2 - \frac{3381}{1000}S_1 + 15\cdot\frac{219}{250} \ge -\frac{50320449}{4000000} = -12.58011225 > U_6.
\]

\paragraph{Exclusion of $k=1$ and non-matching $k=2$ auxiliary classes.}
The short-edge intervals can start at $q=289/400$: by~\eqref{eq:loss-identity}, $E_6(X)\ge f(q_{ij})-14$ for every pair, and $f(289/400)-14>U_6$. Since $f$ decreases on $(0,1]$, this cutoff holds on the entire sublevel $\mathcal{S}_6(U_6)$.
For $k \in \{1,2,3\}$, fixing a canonical true-long matching of size $k$ and exhausting all $2^{15}$ labeled auxiliary masks modulo the stabilizer subgroup of the $k$-matching in $S_6$ yields $18$ auxiliary isomorphism representatives (covering $290$ labeled masks) for $k=1$, $4$ representatives ($14$ labeled masks) for $k=2$, and $1$ matching representative ($1$ mask) for $k=3$. All $18$ single-long classes are excluded by $70$ interval leaves (minimum margin $0.00303501$ above $U_6$), and $3$ of the $4$ two-long classes ($13$ masks) are excluded by $108$ rectangle leaves (minimum margin $0.00009977$ above $U_6$). The single-long leaves use $\lambda=0$ on $61$ leaves and $\lambda=2/25$ on $9$; the double-long leaves use these values on $56$ and $52$ leaves, respectively. The nonzero-multiplier leaves use the stationary virial identity of Proposition~\ref{prop:scale-identity}.
The single remaining $k=2$ class ($1$ mask) has $G_{\mathrm{aux}}$ equal to a $3$-edge matching, which merges with $k=3$ into the matching cube $[3/2, 4B_0]^3$ (Table~\ref{tab:lj6-aux}).

\begin{table}[htbp]
\centering
\caption{Classification of surviving $\mathrm{LJ}_6$ domains and auxiliary graphs ($q_{ij} \ge 3/2$) and their exact exclusion/localization leaf counts (energy margins above $U_6 = -3267/257$ are downward-rounded to six decimal places).}
\label{tab:lj6-aux}
\smallskip
\resizebox{0.98\textwidth}{!}{%
\begin{tabular}{c r r l c r}
\toprule
True-long edges $k$ & Auxiliary classes & Labeled masks & Resolution method & $\lambda$ & Leaves / Min.\ margin \\
\midrule
$0$ & --- (unified) & All on $(0, B_0]^{15}$ & Global Bernstein + rank-3 trace bound & $0$ & $1$ analytical domain ($+0.131950$) \\
$1$ & $18$ & $290$ & Exact rational 1D interval spectral trees & $0,\ 2/25$ & $70$ leaves ($+0.003035$) \\
$2$ & $3$  & $13$  & Exact rational 2D rectangle spectral trees & $0,\ 2/25$ & $108$ leaves ($+0.000099$) \\
$2\text{ or }3$ & $1$ & $2$ & Antipodal matching cube $[3/2, 4B_0]^3$ & $0$ & $3{,}195$ excl.\ ($3{,}123$ spec.\ + $72$ pair) + $40$ loc. \\
\bottomrule
\end{tabular}%
}
\end{table}

\subsection{Explicit Antipodal Gram Witness and the \texorpdfstring{$6{,}469$}{6,469}-Node Matching Tree}
In the surviving topology, $G_{\mathrm{aux}}(X)$ is the matching $\{(0,1),(2,3),(4,5)\}$ with three matching squared distances $q_{2i,2i+1} \in [l_i, h_i] \subseteq [3/2, 4B_0]$ ($i=0,1,2$) and twelve short cross-pair distances $q_{ab} < 3/2$.

Set $K = 2$ ($\lambda = 0$) and choose the rational short contact $z = 9910823276741785/10^{16}$. On $(0, 3/2]$, $f(q) - (Kq^2 + b_s q + c_s) = q^{-6}(q-z)^2 R_6(q)$ with tangent coefficients $b_s, c_s \in \Q$, and all seven Bernstein coefficients of $R_6$ on $[0, 3/2]$ are strictly positive. On each long interval $[l_i, h_i]$, concavity of $f(q) - Kq^2$ provides a secant support $Kq^2 + b_i q + c_i$. Defining $\ell_i := b_i + 2b_s$, $L := \sum_{i=0}^2 \ell_i$, and
\begin{equation}
y_i := -\frac{2}{5K}\left(\ell_i - \frac{L}{8}\right), \qquad \mu_i := \frac{y_i}{2}, \qquad \mathrm{LB} := 12c_s + \sum_{i=0}^2 c_i + \frac{1}{2}\sum_{i=0}^2 \ell_i y_i,
\end{equation}
we construct the antipodal rank-three matrix $T$ with $2\times 2$ blocks $\frac{y_i}{4}\bigl[\begin{smallmatrix} 1 & -1 \\ -1 & 1 \end{smallmatrix}\bigr]$ (eigenvalues $\mu_0, \mu_1, \mu_2 > 0$) and the full $6\times 6$ dual matrix $H$ via Proposition~\ref{prop:gram-identity} with $v = \diag(T - G_0)$ and $\tau = \tr(T - G_0)$. Here $H - T$ vanishes identically on the three antisymmetric pair vectors, and its action on the three pair-center vectors reduces to the $3\times 3$ rational matrix $S_{ij} = 2(H-T)_{2i,2j}$. Bounding $\lambda_{\max}(S) \le b$ via Gershgorin's theorem and verifying $0 \le b < \min_i \mu_i$ ensures that $T$ is the exact top-three truncation of $H$, with spectral gap constant $c = 1 - b/\min_i \mu_i > 0$ and $\mathcal{B}_{\mathrm{spec}}(H) = \mathrm{LB}$.

Branching over $[3/2, 4B_0]^3$ produces a complete binary tree of $6{,}469$ nodes ($3{,}234$ splits and $3{,}235$ leaves): $3{,}195$ leaves are excluded ($3{,}123$ via the spectral bound $\mathrm{LB} > U_6$ and $72$ via the direct pair-loss bound $-12 + \sum_{i=0}^2 f(l_i) > U_6$), leaving $40$ leaves near $(2q_*, 2q_*, 2q_*)$. On these $40$ leaves, both the explicit rank-three truncation $T$ and the exact octahedral Gram matrix $G_{\mathrm{oct}}$ are block-diagonal with $2\times 2$ blocks $\frac{\mu_i}{2}\bigl[\begin{smallmatrix} 1 & -1 \\ -1 & 1 \end{smallmatrix}\bigr]$ and $\frac{q_*}{2}\bigl[\begin{smallmatrix} 1 & -1 \\ -1 & 1 \end{smallmatrix}\bigr]$ ($i=0,1,2$). Since each block $\frac{1}{2}\bigl[\begin{smallmatrix} 1 & -1 \\ -1 & 1 \end{smallmatrix}\bigr]$ has Frobenius norm $1$, enclosing $q_* = (257/264)^{1/3}$ in $[q_{\mathrm{lo}}, q_{\mathrm{hi}}] := [9910823276741784/10^{16},\, 9910823276741787/10^{16}]$ (certified by $q_{\mathrm{lo}}^3 < 257/264 < q_{\mathrm{hi}}^3$) yields directly
\[
\|T - G_{\mathrm{oct}}\|_F^2 = \sum_{i=0}^2 (\mu_i - q_*)^2 \le \delta := \sum_{i=0}^2 \max\bigl(|\mu_i - q_{\mathrm{lo}}|,\, |\mu_i - q_{\mathrm{hi}}|\bigr)^2.
\]
Combining $\|G - T\|_F^2 \le \frac{U_6 - \mathrm{LB}}{2Kc}$ from Proposition~\ref{prop:gram-identity} and Lemma~\ref{lem:spectral-gap} with the squared triangle inequality $\|G - G_{\mathrm{oct}}\|_F^2 \le 2\|G - T\|_F^2 + 2\|T - G_{\mathrm{oct}}\|_F^2$ gives
\[
\|G - G_{\mathrm{oct}}\|_F^2 \le \frac{U_6 - \mathrm{LB}}{Kc} + 2\delta \le 10^{-10}
\]
across all $40$ leaves (maximum $\approx 9.7804 \times 10^{-11}$, certifying $\|G - G_{\mathrm{oct}}\|_F \le 10^{-5}$ and hence $|q_{ij} - q_{ij}^{\mathrm{oct}}| \le 2\times 10^{-5}$ by~\eqref{eq:pair-distance-localization}; by contrast, the generic bound~\eqref{eq:gram-cand-bound} of Corollary~\ref{cor:gram-localization}, which divides $\delta$ by $c$, would exceed $10^{-5}$ on $4$ of the $40$ leaves). Applying the oriented chart reconstruction of Proposition~\ref{prop:chart-newton} (after relabeling $(0,1,2,3,4,5) \mapsto (0,1,2,4,3,5)$ so that the reference tetrahedron is non-coplanar) yields a $12$-dimensional coordinate box of maximum width $\approx 3.6941 \times 10^{-4}$ ($< 3.70 \times 10^{-4}$), which contracts in a single $160$-bit interval Newton step to width $\approx 1.5291 \times 10^{-4}$ ($< 1.53 \times 10^{-4}$), strictly inside the certified octahedral uniqueness box of radius $2^{-13}$. This proves Theorem~\ref{thm:lj6}.

\subsection{Exact Certification and Exclusion of the Competing \texorpdfstring{$C_{2v}$}{C2v} Local Minimum}
In chemical physics, $\mathrm{LJ}_6$ possesses a well-known second local minimum: the $C_{2v}$ capped triangular bipyramid (equivalently, a two-face-capped tetrahedron). To confirm its exact status independently of global branch-and-bound, we parameterize the $C_{2v}$ family by $(\pm a,0,0)$, $(0,b,\pm c)$, and $(0,d,\pm e) \in \R^3$, and certify its unique stationary point in a $5$-variable box of radius $10^{-25}$ and in a full $12$-dimensional reduced gauge box via Krawczyk contraction and positive definite reduced interval Hessians.

\begin{proposition}[Certified $C_{2v}$ Competing Local Minimum]
\label{prop:lj6-c2v}
Modulo rigid Euclidean motions $\mathrm{E}(3)$ and particle permutations $S_6$, the six-particle Lennard--Jones potential possesses a strict local minimum of exact point group $C_{2v}$ (a two-face-capped tetrahedron) whose energy $E_6^{C_{2v}}$ is enclosed within the exact rational bracket
\begin{equation}
-12.302927529580728482747833461462 < E_6^{C_{2v}} < -12.302927529580728482747833461461.
\end{equation}
In particular, $E_6^{C_{2v}} > -3267/257 = E_6^{\min}$ by a certified margin exceeding $0.409134$.
\end{proposition}

\begin{proof}
In the $5$-parameter family $p_{0,1} = (\pm a,0,0)$, $p_{2,3} = (0,b,\pm c)$, and $p_{4,5} = (0,d,\pm e)$, preconditioned Krawczyk contraction and Gershgorin positivity of the $5\times 5$ preconditioned interval Hessian on a rational box $\mathcal{B}_5$ of radius $10^{-25}$ (contraction factor $< 6.07 \times 10^{-19}$) certify a unique critical point $(a^*,b^*,c^*,d^*,e^*) \in \operatorname{int}(\mathcal{B}_5)$ of the restricted energy $\tilde{E}_6(a,b,c,d,e)$. To see that this root is stationary in the full $18$-dimensional configuration space, note that invariance under the two mirror reflections $x \mapsto -x$ and $z \mapsto -z$ forces zero $x$-force on the four $x=0$ particles $(p_2,p_3,p_4,p_5)$, zero $z$-force on the axial pair $(p_0,p_1)$, and equal $y$-forces within each symmetric pair; vanishing of the five partial derivatives $\partial \tilde{E}_6/\partial(a,b,c,d,e) = 0$ eliminates the axial $x$-forces and all $y$- and $z$-forces on $p_2,p_3,p_4,p_5$; and overall translation invariance $\sum_{i=0}^5 \nabla_{p_i} E_6 = 0$ then forces the remaining common $y$-force on the axial pair $(p_0,p_1)$ to vanish. Mapping $\mathcal{B}_5$ into the $12$-dimensional oriented gauge chart $\iota$ of Proposition~\ref{prop:chart-newton} places the entire $5$-parameter box strictly inside a certified $12$-dimensional Krawczyk uniqueness box of radius $10^{-8}$ (contraction factor $< 7.14 \times 10^{-5}$, with positive definite preconditioned $12\times 12$ interval Hessian), identifying the family root with the unique strict local minimizer modulo $\mathrm{E}(3)$ in that neighborhood and yielding the stated energy bracket.

Finally, on $\mathcal{B}_5$, outward interval evaluation shows that no pairwise squared distance straddles $6/5$, and thresholding at $q_{ij} < 6/5$ isolates a fixed $12$-edge contact graph whose automorphism group in $S_6$ has order $4$. Because the six particles in $\mathcal{B}_5$ affinely span $\R^3$, any Euclidean symmetry of the configuration acts faithfully on the vertex set and preserves this $12$-edge contact graph; since the two mirror reflections $x \mapsto -x$ and $z \mapsto -z$ already generate $4$ distinct symmetries, they exhaust the Euclidean point group, proving that the exact point group is $C_{2v}$.
\end{proof}

\section{Complete Proof of Global Optimality for \texorpdfstring{$N=7$}{N=7}}
\label{sec:lj7}

We now present the proof of Theorem~\ref{thm:lj7} for the seven-particle cluster.

\subsection{Exact Algebraic Analysis of the \texorpdfstring{$D_{5h}$}{D5h} Pentagonal Bipyramid}
\label{subsec:lj7-candidate}
Place five equatorial particles at $(R\cos(2\pi k/5), R\sin(2\pi k/5), 0)$ ($k=0,\dots,4$) and two axial apices at $(0,0,\pm H)$ with $R, H > 0$, and define the dimensionless aspect ratio $t := H^2/R^2 > 0$. Dividing the $21$ squared distances by $R^2$ yields five equatorial nearest-neighbor edges of scaled squared length $\alpha = (5-\sqrt{5})/2$, five equatorial diagonals of scaled squared length $\beta = (5+\sqrt{5})/2$, ten apex-equator edges of scaled squared length $1+t$, and one apex-apex edge of scaled squared length $4t$. Because $\alpha^{-1} + \beta^{-1} = 1$ and $(\alpha\beta)^{-1} = 1/5$, Newton's identities yield the exact rational power sums
\begin{equation}
5(\alpha^{-3} + \beta^{-3}) = 2, \qquad 5(\alpha^{-6} + \beta^{-6}) = \frac{18}{25}.
\end{equation}

\begin{table}[htbp]
\centering
\caption{Pairwise squared-distance orbits of the $D_{5h}$ pentagonal bipyramid (divided by $R^2$) and distribution of the $26 = \binom{5}{2}+\binom{5}{3}+\binom{5}{4}+\binom{5}{5}$ surviving labeled true-long subsets on the canonical representative of $G_{\mathrm{aux}} \cong C_5 \sqcup 2K_1$.}
\label{tab:lj7-orbits}
\smallskip
\begin{minipage}[t]{0.54\textwidth}
\centering
{\small \textbf{(a)} $D_{5h}$ candidate distance orbits}\\[4pt]
\small
\setlength{\tabcolsep}{4pt}%
\begin{tabular}{l r c}
\toprule
Orbit type & Pairs & Scaled $q_{ij}/R^2$ \\
\midrule
Equatorial nearest neighbors & $5$  & $\alpha = (5-\sqrt{5})/2$ \\
Equatorial diagonals ($C_5$) & $5$  & $\beta = (5+\sqrt{5})/2$  \\
Apex--equator contacts       & $10$ & $1 + t$                   \\
Apex--apex axis              & $1$  & $4t$                      \\
\bottomrule
\end{tabular}
\end{minipage}
\hfill
\begin{minipage}[t]{0.42\textwidth}
\centering
{\small \textbf{(b)} Surviving $C_5 \sqcup 2K_1$ true-long subsets}\\[4pt]
\small
\setlength{\tabcolsep}{6pt}%
\begin{tabular}{r r}
\toprule
True-long count $k$ & Labeled subsets \\
\midrule
$k = 2$ & $10$ \\
$k = 3$ & $10$ \\
$k = 4$ & $5$  \\
$k = 5$ & $1$  \\
\bottomrule
\end{tabular}
\end{minipage}
\end{table}

Consequently, the $D_{5h}$ family energy is a rational function of $R^6$ and $t$:
\begin{equation}
F(R,t) = A(t)R^{-12} - 2B(t)R^{-6} = A(t)\left(R^{-6} - \frac{B(t)}{A(t)}\right)^2 - \frac{B(t)^2}{A(t)},
\end{equation}
with $A(t)$ and $B(t)$ defined in~\eqref{eq:lj7-AB}. Minimizing over $R > 0$ fixes $R^6 = A(t)/B(t)$ and reduces the energy to $e(t) = -B(t)^2/A(t)$. Differentiating $e(t)$ and clearing common factors yields
\begin{equation}
e'(t) = \frac{3 B(t) P(t)}{25600\,t^7 (1+t)^7 A(t)^2},
\end{equation}
where $P(t) \in \Z[t]$ is the degree-$10$ polynomial in~\eqref{eq:lj7-poly}. An exact rational Sturm sequence shows that $P(t)$ has sign-variation counts $(6, 5, 4)$ at $(-\infty, 0, +\infty)$, proving that $P(t)$ has a unique positive real root $t_* \in (0.36, 0.37)$. With $a(t)$ and $b(t)$ defined in~\eqref{eq:lj7-energy-a} and~\eqref{eq:lj7-energy-b}, respectively, the primitive resultant polynomial $\mathcal{Q}(e) := \operatorname{prim}(\operatorname{Res}_t(P(t), a(t)e + 25b(t)^2)) \in \Z[e]$ has degree $10$, and exact Frobenius irreducibility tests over $\mathbb{F}_{23}$ ($x^{23^{10}} \equiv x \pmod f$ and $\gcd(f, x^{23^d}-x)=1$ for $d \in \{2,5\}$) prove that both $P(t)$ and $\mathcal{Q}(e)$ are irreducible over $\Q$.

By equatorial reflection, vertical mirror reflection through each equatorial vertex, and fivefold rotational symmetry, vanishing of $\partial F/\partial R$ and $\partial F/\partial t$ at $(R_*, t_*)$ implies full $21$-coordinate Cartesian stationarity $\nabla E_7(X^*) = 0$. Furthermore, enclosing the gauge coordinates on a $2^{-192}$ rational grid places $X^*$ near the center of a $15$-dimensional reduced gauge box of radius $2^{-17}$ on which interval Newton contraction (contraction factor $0.038296 < 1$) and positive definiteness of the preconditioned reduced interval Hessian certify local uniqueness and strict local minimality modulo $\mathrm{E}(3)$. Using outward rational bounds on $E_7(X^*)$, we set two rational global thresholds:
\begin{equation}
\label{eq:u0-u1}
E_7^{\min} \le E_7(X^*) < U_7^{(1)} := -16.5053841680 < U_7^{(0)} := -16.5053835.
\end{equation}

\subsection{Distance-Graph Reduction to \texorpdfstring{$C_5 \sqcup 2K_1$}{C5 + 2K1}}
\label{subsec:lj7-graph-reduction}
By Lemma~\ref{lem:attainment} and Remark~\ref{rem:elementary-attainment}, global minimizers of $E_7$ exist in $\mathcal{S}_7(U_7^{(1)}) \subset \mathcal{S}_7(U_7^{(0)})$ and satisfy $\nabla E_7(X) = 0$ (hence $V(X) = 0$). Set the true-long squared-distance threshold to $B_0 = (73/50)^2 = 5329/2500 = 2.1316$ ($r_{\mathrm{long}} = 1.46$) and the auxiliary threshold to $C = (1 + B_0)/2 = 7829/5000 = 1.5658$.

\begin{proposition}[Two-Stage Quotient Cover of the Seven-Particle Distance Graphs]
\label{prop:quotient-cover}
Fix the thresholds $1 < C < B_0$ and the sublevel $\mathcal{S}_7(U_7^{(0)})$. Up to particle relabeling in $S_7$, every configuration $X \in \mathcal{S}_7(U_7^{(0)})$ belongs to the finite inventory constructed by:
\begin{enumerate}[label=(\roman*),leftmargin=2.2em]
    \item partitioning the $2^{21}$ labeled auxiliary edge masks on $V = \{0,\dots,6\}$ into $1{,}044$ unlabeled isomorphism classes under $S_7$, with one canonical representative $G_{\mathrm{aux}}^{\mathrm{rep}}$ per class; and
    \item enumerating on each representative every labeled true-long subset $E_{\mathrm{long}} \subseteq E(G_{\mathrm{aux}}^{\mathrm{rep}})$ satisfying the pruning inequalities of Proposition~\ref{prop:graph-pruning}.
\end{enumerate}
For $2 \le k \le 5$, filtering by the degree and joint-loss budgets leaves $3{,}586$ non-$C_5$ colored cases together with $26 = \binom{5}{2} + \binom{5}{3} + \binom{5}{4} + \binom{5}{5} = 10 + 10 + 5 + 1$ labeled true-long subsets on the canonical representative of $G_{\mathrm{aux}} \cong C_5 \sqcup 2K_1$ (Table~\ref{tab:lj7-orbits}).
\end{proposition}

\begin{proof}
The action of $S_7$ partitions the $2^{21}$ labeled masks into $1{,}044$ isomorphism orbits. For any $X \in \mathcal{S}_7(U_7^{(0)})$, choosing $\sigma \in S_7$ with $G_{\mathrm{aux}}(\sigma X) = G_{\mathrm{aux}}^{\mathrm{rep}}$ preserves all pairwise distances and pruning inequalities, and $E_{\mathrm{long}}(\sigma X) \subseteq E(G_{\mathrm{aux}}^{\mathrm{rep}})$ is one of the enumerated subsets. Since the edge thresholds $q_{ij} \ge C$ and $q_{ij} \ge B_0$ are closed, boundary configurations are covered by adjacent cells.
\end{proof}

Using the six-particle lower bound $E_6(X) \ge -3267/257 > -129/10 =: L_6^{\mathrm{aux}}$ from Theorem~\ref{thm:lj6} (independently cross-checked directly from $L_5$ in the companion Zenodo archive; see Appendix~\ref{app:verification}), Proposition~\ref{prop:graph-pruning} and exact rational spectral verification establish:
\begin{enumerate}[leftmargin=2em]
    \item \textbf{Maximum true-long degree $2$ and edge count $k \le 5$:} Any vertex of true-long degree $3$ forces $E_7 \ge -12.9 - 3 + 3v(1.46) \approx -16.487509 > U_7^{(0)}$, and $6$ true-long edges force $E_7 \ge -15 + 6v(1.46) \approx -16.175018 > U_7^{(0)}$. Moreover, since $k \le 5$ and maximum degree is at most $2$, for every true-long edge $\{i,j\}$ at least $3$ of the remaining $5$ vertices are connected to both $i$ and $j$ by non-long edges ($< r_{\mathrm{long}}$), so $q_{ij} \in [B_0, 4B_0]$ by~\eqref{eq:root-box-enclosure}.
    \item \textbf{Stationary spectral exclusion of $k=0$ (zero true-long edges):} On $(0, B_0]$, for a stationary virial multiplier $\lambda \in \Q$ ($V(X)=0$) and two-piece quadratic weight $\gamma \in \Q_{>0}$, we certify a two-piece rational Bernstein minorant of the form $f_\lambda(q) \ge a q^2 + b q + c + \gamma\min\{(q-1)^2, (q-B_0)^2\}$ with common leading quadratic coefficient $K = a + \gamma$ across the two regimes $q \le C$ and $q \ge C$. When all $21$ squared distances satisfy $q_{ij} \le B_0$, taking $(\lambda, \gamma) = (1/2, 2)$ splits $(0, B_0]^{21}$ across the $1{,}044$ unlabeled auxiliary graph representatives (verified to cover all $2^{21}$ labeled masks). Exact rational $LDL^\top$ spectral verification across all $1{,}044$ representatives proves $E_7(X) > -16.14135379 > U_7^{(0)}$ (weakest lower bound $\approx -16.141353783$).
    \item \textbf{Stationary spectral exclusion of $k=1$ (single true-long edge):} When $k=1$, the joint loss inequality~\eqref{eq:joint-edge-budget} gives $\ell(B_0) + (m-1)\ell(C) \le 21 + U_7^{(0)}$, so $G_{\mathrm{aux}}(X)$ has at most $m \le \lfloor 1 + (21 + U_7^{(0)} - \ell(B_0))/\ell(C)\rfloor = 7$ edges (since $m \ge 8$ forces $E_7 \ge -21 + \ell(B_0) + 7\ell(C) \approx -16.36771157 > U_7^{(0)}$). Selecting a distinguished true-long edge in each of the $145$ non-empty auxiliary isomorphism representatives with $1 \le m \le 7$ edges yields $866$ marked-edge cases; dividing the single true-long squared distance $[B_0, 4B_0]$ into $24$ slabs produces $24 \times 866 = 20{,}784$ slab-graph records. Exact verification of all $20{,}784$ rational spectral records (using two certified two-piece minorant profiles on the $20$ non-long pairs: $(\lambda, \gamma) = (1/2, 2)$ on $20{,}781$ records and $(\lambda, \gamma) = (1/4, 159/100)$ on $3$ records) proves $E_7(X) > -16.45048964 > U_7^{(0)}$ (weakest lower bound $\approx -16.450489635$).
    \item \textbf{Complete exclusion of all $3{,}586$ non-$C_5$ colored representative cases ($2 \le k \le 5$):} Filtering all $1{,}044$ unlabeled auxiliary representatives and all labeled true-long subsets $E_{\mathrm{long}} \subseteq E(G_{\mathrm{aux}}^{\mathrm{rep}})$ of size $2 \le k \le 5$ with maximum degree at most $2$ by the loss inequality~\eqref{eq:joint-edge-budget} yields $3{,}586$ non-$C_5$ colored cases and $26 = \binom{5}{2} + \binom{5}{3} + \binom{5}{4} + \binom{5}{5} = 10 + 10 + 5 + 1$ labeled true-long subsets on the canonical representative of $G_{\mathrm{aux}} \cong C_5 \sqcup 2K_1$ (Table~\ref{tab:lj7-orbits}). Independent from-root rational spectral replay closes all $3{,}586$ non-$C_5$ cases in $129{,}506$ events ($62{,}960$ splits and $66{,}546$ exclusion leaves), proving $E_7(X) > -16.50534781 > U_7^{(0)}$ (weakest lower bound $\approx -16.505347804$).
\end{enumerate}

Consequently, up to relabeling in $S_7$, every stationary configuration with $E_7(X) \le U_7^{(0)}$ must have auxiliary graph $G_{\mathrm{aux}}(X) \cong C_5 \sqcup 2K_1$: five vertices form an auxiliary $5$-cycle of edges with $q_{ij} \ge C$, and the remaining two vertices (the apices $0,1$) have no incident auxiliary edges ($q_{ab} < C$ for all $16$ non-auxiliary pairs).

\subsection{Closing the Auxiliary \texorpdfstring{$C_5$}{C5} Domain: Root Geometry, Non-Strong, and Strong Regimes}
\label{subsec:lj7-c5-closure}
Let $q_{e_1},\dots,q_{e_5} \ge C$ denote the five squared distances along the auxiliary $5$-cycle $C_5$, and let $a,b \in \R^3$ denote the two apex positions ($0,1$) with midpoint $m = (a+b)/2$ and squared separation $q_{01} = \|a-b\|_2^2$. To enclose $q_{e_1},\dots,q_{e_5}$ in a finite root box without assuming coplanarity of the $5$-cycle or collinearity of the apices with its center, observe that if every cycle vertex $x$ satisfies $\|x-a\|_2^2 \le s^2$ and $\|x-b\|_2^2 \le s^2$ for some upper bound $s > 0$, the parallelogram identity gives
\begin{equation}
\label{eq:midpoint-parallelogram}
\|x - m\|_2^2 = \frac{\|x-a\|_2^2 + \|x-b\|_2^2}{2} - \frac{\|a-b\|_2^2}{4} \le s^2 - \frac{q_{01}}{4}.
\end{equation}
Hence, by the triangle inequality through $m$, any two cycle vertices $x,y$ satisfy
\begin{equation}
\label{eq:cycle-chord-bound}
q_{xy} = \|x - y\|_2^2 \le \bigl(\|x-m\|_2 + \|y-m\|_2\bigr)^2 \le 4s^2 - q_{01}.
\end{equation}
We now partition the $C_5 \sqcup 2K_1$ domain into two complementary regimes that overlap on the boundary $\min_{1 \le m \le 5} q_{e_m} = B_0$:

\paragraph{1. The non-strong $C_5$ regime ($\min_{1 \le m \le 5} q_{e_m} \le B_0$).}
Since all five cycle edges satisfy $q_{e_m} \ge C = 7829/5000$, any non-cycle pair $(u,v)$ obeys $\ell(q_{uv}) \le 21 + U_7^{(0)} - 5\ell(C)$, which enforces $r_{01} = \sqrt{q_{01}} > 217/250 = 0.868$. The ten apex--equator edges are non-auxiliary, so every cycle vertex $x$ obeys $\|x-a\|_2^2 < C$ and $\|x-b\|_2^2 < C$; this is precisely the hypothesis of~\eqref{eq:midpoint-parallelogram}--\eqref{eq:cycle-chord-bound} with $s^2 = C$, which therefore bounds every cycle squared distance above by $4C - q_{01} < 4C - (217/250)^2 = 344361/62500 = 5.509776$. If at least one cycle edge has squared length at most $B_0$, dihedral symmetry of $C_5$ allows us to fix $q_{e_1} = \min_{m} q_{e_m} \in [C, B_0]$ while retaining $q_{e_m} \in [q_{e_1}, 344361/62500]$ for $m=2,\dots,5$, with the apex-apex distance $q_{01}$ enclosed as a derived interval from~\eqref{eq:cycle-chord-bound} and coupled via six selected-pair interval multiplier constraints ($q_{01}, q_{e_1},\dots,q_{e_5}$). Five-dimensional branching over $(q_{e_1},\dots,q_{e_5})$ closes the entire non-strong regime in $1{,}767$ events ($883$ splits and $884$ exact spectral exclusion leaves), with weakest certified lower bound $-16.505092607\dots$, exceeding $U_7^{(1)}$ by $2.9156 \times 10^{-4}$ and $U_7^{(0)}$ by $2.9089 \times 10^{-4}$.

\paragraph{2. The strong five-long $C_5$ regime ($q_{e_m} \ge B_0$ for all $m=1,\dots,5$).}
When all five auxiliary cycle edges are true-long edges ($q_{e_m} \ge B_0$), the loss budget~\eqref{eq:loss-identity} bounds the loss of each of the $16$ non-cycle pairs (the apex--apex edge, the ten apex--equator edges, and the five equatorial nearest-neighbor edges) by $\ell(q_{ab}) \le 21 + U_7^{(0)} - 5\ell(B_0) \approx 0.473798348$, which confines every short distance to $0.916 < 0.916411 < r_{ab} < 1.214463 < 1.215$. In particular, every cycle vertex lies within distance $1.215$ of each apex, supplying the hypothesis of~\eqref{eq:cycle-chord-bound} with $s = 1.215$; together with the apex-apex enclosure $q_{01} > (0.916)^2$, this gives $q_{e_m} < 4(1.215)^2 - (0.916)^2 = 1266461/250000 = 5.065844$, enclosing the five cycle squared distances in the root hypercube
\begin{equation}
\label{eq:lj7-strong-cube}
(q_{e_1},\dots,q_{e_5}) \in \left[\frac{5329}{2500},\; \frac{1266461}{250000}\right]^5 = [2.1316,\; 5.065844]^5.
\end{equation}
To certify~\eqref{eq:lj7-strong-cube} with $K = 1$ and $\lambda = 0$, the interval-dual relaxation combines concavity secants on the five long cycle edges $e_1,\dots,e_5$ with three orbit-specific exact rational Bernstein tangent quadratic supports on the $16$ short pairs (enclosed in $[(229/250)^2, (243/200)^2]$): \texttt{apex\_pair} on the apex--apex edge $\{0,1\}$, \texttt{apex\_ring} on the ten apex--equator edges, and \texttt{ring\_short} on the five equatorial nearest-neighbor edges, together with nonnegative endpoint multipliers $(\alpha,\beta)$ on the six selected pairs $\{0,1\}, e_1,\dots,e_5$. The tree checker reconstructs every leaf box recursively from the root hypercube~\eqref{eq:lj7-strong-cube} across all $153{,}333$ schedule events ($76{,}666$ binary splits and $76{,}667$ leaves):
\begin{itemize}[leftmargin=2em]
    \item \textbf{$72{,}898$ exclusion leaves:} Certified strictly above $U_7^{(1)}$ via the above $K=1, \lambda=0$ supports, endpoint multipliers, and exact rational $LDL^\top$ inertia checks on $H$.
    \item \textbf{$3{,}769$ localized leaves:} The checker verifies by exact rational interval comparison that each of the $3{,}769$ reconstructed leaf boxes is contained in a single shared $5$-dimensional box \texttt{candidate\_hull} ($[184760386426845/2^{46},\, 2956210291516585/2^{50}]^5$), and certifies one shared spectral-localization witness on \texttt{candidate\_hull}. Specifically, $LDL^\top$ inertia on \texttt{candidate\_hull} certifies $\lambda_3(H) \ge a = 13/25$ and $\lambda_4(H) \le b = 11/25$ ($c = 1 - b/a = 2/13$); substituting $(U_7^{(1)}, \mathrm{LB}, K=1, c=2/13, \delta)$ into~\eqref{eq:gram-cand-bound} of Corollary~\ref{cor:gram-localization} yields the uniform Frobenius radius
    \begin{equation}
    \label{eq:lj7-strong-localization}
    \|G - G_{\mathrm{cand}}\|_F \le \frac{6555377689298343137}{75557863725914323419136} < 8.675971 \times 10^{-5} < 10^{-4}.
    \end{equation}
\end{itemize}
Using~\eqref{eq:lj7-strong-localization} on \texttt{candidate\_hull} and applying Proposition~\ref{prop:chart-newton} in the apex-axis gauge (with particles $0,1$ defining the positive axis) recovers a $15$-dimensional reduced Cartesian box of maximum width $\approx 4.5040 \times 10^{-3}$ ($< 4.51 \times 10^{-3}$). Four preconditioned raw-Hessian interval Newton contractions $\mathcal{N}_k$ (using direct pair-correlated enclosures $\mathcal{M}_{Y_{k-1}}(\mathcal{X}_{k-1}) = [Y_{k-1}\nabla^2 F(\mathcal{X}_{k-1})]$ on the full collision-free coordinate boxes) map this box strictly inside the $2^{-17}$-radius uniqueness box $\mathcal{X}_{\mathrm{loc}}$ of the exact pentagonal bipyramid $X^*$. Since every global minimizer of $E_7$ must enter this unique-critical-point box, Theorem~\ref{thm:lj7} is proved.

\section{Propagation to \texorpdfstring{$8 \le N \le 150$}{8 <= N <= 150} and the \texorpdfstring{$N=8$}{N=8} Frontier}
\label{sec:propagation-lj8}

\subsection{Subset Energy Envelopes and Certified Upper Bounds}
\label{subsec:subset-envelopes}
For any $N \ge M \ge 2$, summing the certified $M$-particle lower bound $E_M(X_S) \ge L_M$ over all $\binom{N}{M}$ subsets $S \subset \{1,\dots,N\}$ double-counts each pair $\{i,j\}$ exactly $\binom{N-2}{M-2}$ times, yielding the classical subset scaling bound
\begin{equation}
\label{eq:subset-scaling}
E_N(X) = \frac{1}{\binom{N-2}{M-2}} \sum_{|S|=M} E_M(X_S) \ge \frac{N(N-1)}{M(M-1)}\,L_M.
\end{equation}
Substituting our newly certified constants $L_6 = -3267/257$ (Theorem~\ref{thm:lj6}) and $L_7 = -16.5053841681$ (Theorem~\ref{thm:lj7}), the seven-particle per-pair coefficient $L_7/42 \approx -0.3929853373$ is strictly greater than the six-particle coefficient $L_6/30 \approx -0.4237354086$. Combining~\eqref{eq:subset-scaling} for $M=7$ with Yuhjtman's analytical stability bound $E_N(X) \ge -14.316\,N$~\cite[Theorem~5.1]{yuhjtman2015poignant} (which holds for every finite configuration in $\Omega_N$ under our potential normalization~\eqref{eq:lj-energy}) establishes Theorem~\ref{thm:propagation}(i), strictly improving the combined six-particle and stability envelope for all $8 \le N \le 37$ (for $N \ge 38$, $-14.316\,N$ dominates both subset bounds).

To obtain certified rational feasible upper bounds $U_N \ge E_N^{\min}$ for $8 \le N \le 150$ without relying on floating-point energy tables, we convert each putative coordinate set in the Cambridge Cluster Database~\cite{cambridgeCCD} from $\sigma$-units to equilibrium-distance units using the rational multiplier $\kappa = 8908987181403393 / 10^{16} \approx 2^{-1/6}$, round coordinates to the $10^{-12}$ rational grid, bound all $\binom{N}{2}$ rational pair energies upward on a $10^{-30}$ integer grid, and round the sum upward to $9$ decimal places. Because $U_N$ is evaluated directly from an explicit rational configuration in $\Omega_N$, it provides a constructively certified upper bound on $E_N^{\min}$.

\subsection{Arrangement-Free Layer Diameter Bounds and Kernel Correction}
\label{subsec:layer-bounds}
Let $X \in \Omega_N$ be a global minimizer of $E_N$ with diameter $d = \max_{i<j}\|x_i - x_j\|_2$. Rotating the diameter pair onto the $x$-axis with left endpoint at $x=0$ projects all particles into $[0,d]$. Partition $[0,\infty)$ into disjoint half-open unit layers $\mathcal{L}_m := [m, m+1)$ for $m = 0,1,\dots, D-1$, where $D := \lfloor d \rfloor + 1 > d$. Every internal layer $\mathcal{L}_m$ ($1 \le m \le D-2$) must be non-empty: otherwise, letting $g > 1$ denote the minimum $x$-separation between particles in $\bigcup_{r \le m-1}\mathcal{L}_r$ and particles in $\bigcup_{r \ge m+1}\mathcal{L}_r$, translating all particles in $\bigcup_{r \ge m+1} \mathcal{L}_r$ leftward by any $\varepsilon \in \bigl(0, \min(1, g-1)\bigr)$ strictly decreases every cross-gap distance while keeping all cross-gap distances strictly above $1$, strictly lowering $E_N(X)$. Since $\mathcal{L}_0$ and $\mathcal{L}_{D-1}$ contain the diameter endpoints, $n_m := |\mathcal{L}_m \cap X| \ge 1$ for all $0 \le m \le D-1$, so $D \le N$.

For two particles in layers $\mathcal{L}_a$ and $\mathcal{L}_b$ with index separation $j = |b - a| \ge 1$, their $x$-separation is strictly greater than $j - 1$. Since $v(r)$ is strictly increasing on $[1,\infty)$ and $\ge -1$ everywhere, their pair interaction is bounded below by the inter-layer kernel
\begin{equation}
\label{eq:layer-kernel}
w(j) := \begin{cases}
-1, & j \in \{1, 2\}, \\
v(j-1) = (j-1)^{-12} - 2(j-1)^{-6}, & j \ge 3.
\end{cases}
\end{equation}
\begin{remark}[Inter-Layer Kernel Indexing Convention]
\label{rem:kernel-correction}
In Kuznetsov and Sahinidis~\cite{kuznetsov2025new}, Table~1 and Algorithm~1 specify the inter-layer kernel for layer separation $j \ge 3$ as $v(j-1)$ and record the last excluded layer count, whereas the standalone C utility \texttt{distance\_ub.c} in the companion repository initializes the kernel at $v(j)$ (which is not a valid lower bound for separation $j$, since points at $x = 0.99 \in \mathcal{L}_0$ and $x = 3.00 \in \mathcal{L}_3$ have distance $2.01$ with $v(2.01) \approx -0.03010 < v(3) \approx -0.00274$) and outputs the first unexcluded layer count. Throughout our comparisons, we use the valid kernel~\eqref{eq:layer-kernel} and the half-open layer convention $d < D$ uniformly.
\end{remark}

To bound the internal and cross-layer energy for a trial layer count $D \in \{2,\dots,N\}$ without assuming any conjecture on the optimal ordering of layer occupancies, we use the arrangement-free relaxation of Lecomte~\cite[Section~5, Lemma~2]{lecomte2026certified}. Write $n_m = 1 + x_m$ with surplus $x_m \ge 0$ and $\sum_{m=0}^{D-1} x_m = N - D =: e$. Define the single-particle background energy and minimum central-layer field
\begin{equation}
\label{eq:layer-background-defs}
B_D := \sum_{j=1}^{D-1} (D - j)\,w(j), \qquad
c_D := \sum_{j=1}^{\lfloor(D-1)/2\rfloor} w(j) + \sum_{j=1}^{\lceil(D-1)/2\rceil} w(j),
\end{equation}
and bound the internal energy excess of a layer with $1 + x_a$ particles by $f_{\mathrm{sub}}(1) := 0$ and $f_{\mathrm{sub}}(n) := \max_{M \in \{2,5,6,7\},\, M \le n} \frac{n(n-1)}{M(M-1)}L_M$ for $n \ge 2$. The propagated bound uses $L_2=-1$, $L_5=-71123847/7812500$, $L_6=-3267/257$, and $L_7=-16.5053841681$.
For any positive integer partition $x = (x_1 \ge \dots \ge x_k \ge 1)$ of $e$ ($x \vdash e$; for $e=0$ the empty sum is $0$), let $\pi_1 \ge \dots \ge \pi_{\binom{k}{2}}$ be the nonincreasing ordering of the products $\{x_a x_b : 1 \le a < b \le k\}$, and let $\gamma_1 \le \dots \le \gamma_{\binom{k}{2}}$ be the nondecreasing ordering of the contiguous-site distance multiset $\{1^{\times(k-1)}, 2^{\times(k-2)}, \dots, (k-1)^{\times 1}\}$. Because any $k$ distinct occupied surplus layers $0 \le p_1 < \dots < p_k \le D-1$ satisfy $p_b - p_a \ge b - a$, the sorted layer separations dominate $\gamma_r$ termwise; since $-w(j) \ge 0$ is nonincreasing, the rearrangement inequality bounds the cross-surplus attraction above by $\sum_{r=1}^{\binom{k}{2}} \pi_r(-w(\gamma_r))$, while each layer's background field $\beta_m := \sum_{n \ne m} w(|n-m|) \ge c_D$ yields $\sum_{m=0}^{D-1} x_m \beta_m \ge e\,c_D$, giving the arrangement-free lower bound
\begin{equation}
\label{eq:arrangement-free-layer}
E_N(X) \ge B_D + e\,c_D + \min_{x \vdash e}\left[\sum_{a=1}^k f_{\mathrm{sub}}(1 + x_a) + \sum_{r=1}^{\binom{k}{2}} \pi_r w(\gamma_r)\right].
\end{equation}
(Note that including partitions with $k > D$ or bounding $\sum_{m=0}^{D-1} x_m \beta_m$ by $e\,c_D$ simply relaxes the lower bound, and at the single-cluster boundary partition $x = (e)$ the sorted background sum coincides with $e\,c_D$.) If every $D \in \{D_0,\dots,N\}$ yields a lower bound strictly above $U_N$, then $\lfloor d \rfloor + 1 \le D_0 - 1$, proving $d < D_0 - 1$. Across $8 \le N \le 150$,~\eqref{eq:arrangement-free-layer} with $L_7$ improves $90$ cluster sizes over the pair-layer baseline and $24$ cluster sizes over the $L_6$-layer baseline (note that Lecomte's count of $92$ improved sizes in~\cite{lecomte2026certified} refers to a different range $5 \le N \le 200$ and baseline convention).

\begin{remark}[Sorted Background Fields]
\label{rem:sorted-layer-background}
Instead of bounding every layer background field $\beta_m$ by its minimum $c_D$, one may sort the $D$ values $\beta_m$ as $\beta_{(1)} \le \cdots \le \beta_{(D)}$. For a positive surplus partition $x_1 \ge \cdots \ge x_k > 0$ with $k \le D$, the rearrangement inequality gives $\sum_{a=1}^k x_a\beta_{p_a} \ge \sum_{a=1}^k x_a\beta_{(a)} \ge e\,c_D$ for every placement in distinct layers $p_a$, yielding the strengthened lower bound
\begin{equation}
\label{eq:sorted-background-layer}
E_N(X) \ge B_D +
\min_{\substack{x\vdash e\\ k\le D}}
\left[\sum_{a=1}^k f_{\mathrm{sub}}(1+x_a)
+\sum_{a=1}^k x_a\beta_{(a)}
+\sum_{r=1}^{\binom{k}{2}}\pi_r w(\gamma_r)\right],
\end{equation}
with the empty partition retained when $e=0$. Across all $143$ cluster sizes $8\le N\le150$,~\eqref{eq:sorted-background-layer} raises the $L_7$ energy lower bound at the first unexcluded layer count for $135$ sizes (by up to $0.0282643$ at $N=11$), while leaving every integer diameter bound unchanged (the smallest remaining gap below $U_N$ is $\approx 0.150430$ at $N=123$).
\end{remark}

\begin{figure}[htbp]
\centering
\includegraphics[width=0.95\textwidth]{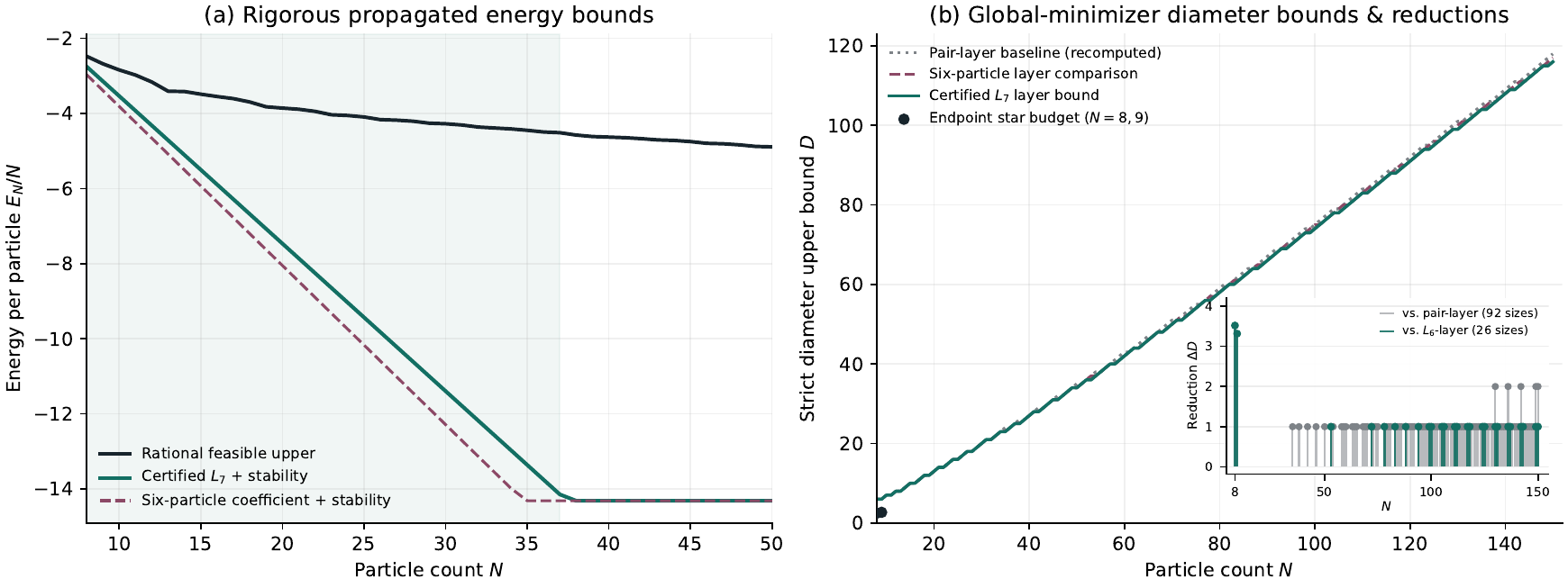}
\caption{Left: Certified per-particle energy bounds ($E_N/N$) across $8 \le N \le 50$, showing the strict improvement of the $\mathrm{LJ}_7$ subset envelope over the $\mathrm{LJ}_6$ subset envelope for $8 \le N \le 37$ before joining Yuhjtman's stability bound. Right: Certified global-minimizer diameter upper bounds across $8 \le N \le 150$, with an inset displaying the exact diameter reductions $\Delta D$ relative to the pair-layer and $L_6$-layer baselines across $8 \le N \le 150$; arrangement-free layer enumeration improves $90$ cluster sizes over the pair-layer baseline and $24$ cluster sizes over the $\mathrm{LJ}_6$-layer baseline ($92$ and $26$ sizes in total when combined with the continuous $N=8,9$ endpoint star bounds).}
\label{fig:bounds}
\end{figure}

For $N=8$ and $N=9$, we obtain even sharper continuous diameter bounds valid across the entire sublevel $\mathcal{S}_N(U_N)$ via two-endpoint star removal. Suppose $X \in \mathcal{S}_N(U_N)$ has diameter $d = \|x_a - x_b\|_2 \ge \rho > 2$. For every $i \notin \{a,b\}$, $r_{ai} + r_{bi} \ge d \ge \rho$, so $v(r_{ai}) + v(r_{bi}) \ge h(\rho) := \inf\{v(r_1)+v(r_2) : r_1+r_2 \ge \rho,\, r_1,r_2 > 0\} = \min_{r_1 \in [1,\rho-1]} g_\rho(r_1)$, where $g_\rho(r_1) := v(r_1)+v(\rho-r_1)$. To see the equality, raise any $r_j<1$ to $1$; this preserves $r_1+r_2\ge\rho$ and cannot increase the energy sum. As $\rho>2$, the distances can then be reduced to sum to $\rho$ while both remain at least $1$, again without increasing the sum because $v$ is increasing on $[1,\infty)$.
The function $h(\rho)$ is nondecreasing in $\rho > 2$; partitioning $[1,\rho-1]$ into $1{,}024$ equal rational subintervals $[l_k, u_k]$ ($0 \le k < 1024$), monotonicity of $v$ on $[1,\infty)$ gives $g_\rho(r_1) \ge v(l_k) + v(\rho - u_k)$ on $[l_k, u_k]$, and hence $h(\rho) \ge \min_{0 \le k < 1024}\bigl(v(l_k) + v(\rho - u_k)\bigr)$. Combining the two-point deletion bound $E_N(X) = E_{N\setminus\{a,b\}} + v(d) + \sum_{i \notin \{a,b\}}(v(r_{ai})+v(r_{bi}))$ with the averaged single-endpoint identity
\begin{equation}
\label{eq:endpoint-identity}
E_N(X) = \frac{1}{2}\Bigl(E_{N\setminus\{a\}}(X_{N\setminus\{a\}}) + E_{N\setminus\{b\}}(X_{N\setminus\{b\}})\Bigr) + v(d) + \frac{1}{2}\sum_{i \notin \{a,b\}}\bigl(v(r_{ai}) + v(r_{bi})\bigr)
\end{equation}
yields, for every configuration with $d \ge \rho > 2$,
\begin{equation}
\label{eq:endpoint-star}
E_N(X) \ge \max\left\{\, L_{N-2} + v(\rho) + (N-2)h(\rho),\; L_{N-1} + v(\rho) + \frac{N-2}{2}h(\rho)\,\right\}.
\end{equation}
At $\rho = 24827/10000 = 2.4827$ for $N=8$ and $\rho = 13429/5000 = 2.6858$ for $N=9$, the right-hand side of~\eqref{eq:endpoint-star} exceeds $U_8$ and $U_9$ by certified margins $1.3109 \times 10^{-4}$ and $1.1165 \times 10^{-4}$, certifying $D_8 < 2.4827$ and $D_9 < 2.6858$ on all of $\mathcal{S}_8(U_8)$ and $\mathcal{S}_9(U_9)$. Table~\ref{tab:propagation} and Figure~\ref{fig:bounds} summarize the certified bounds across $8 \le N \le 150$.

\begin{table}[htbp]
\centering
\caption{Selected certified energy and diameter bounds for representative cluster sizes in $8 \le N \le 150$. The first data column is the certified rational feasible upper bound $U_N$; energy lower bounds combine the subset envelope with Yuhjtman's stability bound (printed downward-rounded to $9$ decimals), and diameter entries are exact. Bold entries mark strict improvements of the $L_7$ bounds over the corresponding $L_6$ bounds (energy lower bounds for $8 \le N \le 37$, layer diameter bounds for $24$ sizes, and final diameter bounds for $26$ sizes when combined with the continuous $N=8,9$ endpoint star bounds).}
\label{tab:propagation}
\smallskip
\resizebox{0.99\textwidth}{!}{%
\begin{tabular}{r r r r r r r r}
\toprule
$N$ & Certified $U_N$ & $L_6$/stab.\ lower & $L_7$/stab.\ lower & Pair-layer $D$ & $L_6$-layer $D$ & $L_7$-layer $D$ & Final $D$ \\
\midrule
  8 & $-19.821489187$ & $-23.729182880$ & $\mathbf{-22.007178891}$ & $6$ & $6$ & $6$ & $\mathbf{2.4827}$ \\
  9 & $-24.113360433$ & $-30.508949417$ & $\mathbf{-28.294944289}$ & $6$ & $6$ & $6$ & $\mathbf{2.6858}$ \\
 35 & $-155.756643135$ & $-501.060000000$ & $\mathbf{-467.652551430}$ & $24$ & $23$ & $23$ & $23$ \\
 38 & $-173.928426511$ & $-544.008000000$ & $-544.008000000$ & $26$ & $25$ & $25$ & $25$ \\
100 & $-557.039819251$ & $-1431.600000000$ & $-1431.600000000$ & $75$ & $75$ & $\mathbf{74}$ & $\mathbf{74}$ \\
150 & $-893.310257879$ & $-2147.400000000$ & $-2147.400000000$ & $118$ & $117$ & $\mathbf{116}$ & $\mathbf{116}$ \\
\bottomrule
\end{tabular}%
}
\end{table}

\subsection{The \texorpdfstring{$N=8$}{N=8} Frontier: Local Certificates, Graph Reduction, and Obstructions}
\label{subsec:lj8-frontier}
Finally, we apply our exact certification tools to investigate the $N=8$ frontier:
\begin{enumerate}[leftmargin=2em]
    \item \textbf{Certified $C_s$ ground-state candidate and a proposed $D_{2d}$ branch:} In an $18$-dimensional reduced rigid-motion gauge box of radius $10^{-6}$, interval Cholesky--Gershgorin Hessian checks and Krawczyk contraction (norm $< 0.01424$) certify a unique strict local minimizer (modulo $\mathrm{E}(3)$) of exact point group $C_s$ (a capped pentagonal bipyramid) with energy enclosed in width $9.31 \times 10^{-46}$ around
    \[
    E_8^{C_s} \approx -19.821489192154760636229140601244399040831505629,
    \]
    providing the rational feasible upper bound $U_8^{\mathrm{loc}} = -1982148919215/10^{11}$. This is the tighter cutoff used for the graph analysis below; Theorem~\ref{thm:propagation}(ii) and Table~\ref{tab:propagation} use the CCD witness bound $U_8=-19.821489187$. A separate $18$-variable local certificate (contraction norm $< 0.02011$, with $8$ certified distance-compatible permutations) encloses a strict local minimum modulo $\mathrm{E}(3)$ on the proposed $D_{2d}$ branch near $E_8^{D_{2d}} \approx -19.76529785090646$, strictly above $U_8^{\mathrm{loc}}$ by $0.056191$; the point group $D_{2d}$ is consistent with these eight permutations, whereas exact point-group invariance is certified for the $C_s$ ground-state candidate.
    \item \textbf{Distance-graph reduction using $L_6$ and $L_7$:} Applying Proposition~\ref{prop:graph-pruning}(i)--(ii) on the sublevel $E_8(X) \le U_8^{\mathrm{loc}}$ with $L_7 = -16.5053841681$ and $L_6 = -3267/257$ restricts $G_{\mathrm{long}}(X)$ to maximum degree $4$ and at most $9$ edges at $r_{\mathrm{long}} = 1.5$ ($10$ edges at $r_{\mathrm{long}}=1.46$). Combining single-vertex and two-vertex deletion budgets reduces the number of unlabeled true-long graphs on the general sublevel $E_8(X) \le U_8^{\mathrm{loc}}$ to $1{,}074$ at $r_{\mathrm{long}}=1.46$, $681$ at $r_{\mathrm{long}}=1.5$, and $514$ at $r_{\mathrm{long}}=1.6$ (Table~\ref{tab:lj8-graphs}).
    \item \textbf{Stationary zero-long exclusion and the one-long spectral bottleneck:} For stationary configurations ($\nabla E_8(X) = 0$, so $V(X) = 0$), at $r_{\mathrm{long}} = 73/50 = 1.46$ ($B_0 = 5329/2500$), taking $\lambda = 2$ in Proposition~\ref{prop:scale-identity} and using the $8$-point rank-three Gram inequality $\sum_{i<j} q_{ij}^2 \ge \frac{1}{24}(\sum_{i<j} q_{ij})^2$ excludes the zero-long graph $G_0$ globally with stationary lower bound $E_8 \ge -17.35610173 > U_8^{\mathrm{loc}}$ (which reduces the stationary graph inventory at $r_{\mathrm{long}}=1.46$---and hence also on the smaller domain $r_{\mathrm{long}}=1.40$---from $1{,}074$ to $1{,}073$ unlabeled graphs and from $19{,}044{,}225$ to $19{,}044{,}224$ labeled graphs). However, on the single-long interval $[B_0, 4B_0]$ partitioned into $48$ slabs, uncolored rank-three Gram relaxations exclude $24$ slabs while leaving the central interval $q_{01} \in [2.39805, 5.59545]$ open (weakest lower bound $-24.8942$), and the independent auxiliary-graph inventory at $C = 1.5658$ (under single- and two-vertex deletion budgets prior to two-color pairing) already comprises $6{,}850$ unlabeled graphs ($152{,}688{,}956$ labeled graphs). Within our spectral-Gram framework, advancing to a complete $N=8$ proof therefore motivates a strengthened auxiliary-colored or subcluster-coupled relaxation before large-scale tree search.
\end{enumerate}

\begin{table}[htbp]
\centering
\caption{Exact inventory of necessary true-long distance graphs for $\mathrm{LJ}_8$ on the general sublevel $E_8(X) \le U_8^{\mathrm{loc}} = -19.82148919215$ (single-vertex and two-vertex deletion budgets) and on the stationary sublevel ($\nabla E_8(X) = 0$, where the $\lambda = 2$ rank-three Gram bound additionally excludes the empty graph $G_0$ for $r_{\mathrm{long}} \le 1.46$).}
\label{tab:lj8-graphs}
\smallskip
\resizebox{0.99\textwidth}{!}{%
\begin{tabular}{r r r r r r r r}
\toprule
Threshold $r_{\mathrm{long}}$ & Max deg. & Max edges & Unlabeled (1-vtx) & Labeled (1-vtx) & Unlabeled (2-vtx) & Labeled (2-vtx) & Stationary (Unlab.\ / Lab.) \\
\midrule
$1.40$ & $4$ & $10$ & $1{,}134$ & $19{,}806{,}525$ & $1{,}074$ & $19{,}044{,}225$ & $1{,}073^*$ / $19{,}044{,}224^*$ \\
$1.46$ & $4$ & $10$ & $1{,}134$ & $19{,}806{,}525$ & $1{,}074$ & $19{,}044{,}225$ & $1{,}073^*$ / $19{,}044{,}224^*$ \\
$1.50$ & $4$ & $9$  & $696$     & $10{,}347{,}299$ & $681$     & $10{,}246{,}499$ & $681$ / $10{,}246{,}499$ \\
$1.60$ & $4$ & $9$  & $696$     & $10{,}347{,}299$ & $514$     & $7{,}390{,}499$  & $514$ / $7{,}390{,}499$ \\
$1.80$ & $3$ & $8$  & $236$     & $2{,}604{,}879$  & $236$     & $2{,}604{,}879$  & $236$ / $2{,}604{,}879$ \\
\bottomrule
\end{tabular}%
}
\\[3pt]
{\footnotesize $^*$Excludes the empty graph $G_0$ via the stationary ($\lambda = 2$) rank-three Gram inequality certified at $r_{\mathrm{long}} = 1.46$ (which contains $r_{\mathrm{long}} = 1.40$).}
\end{table}

\section{Conclusion}
\label{sec:conclusion}

By reformulating the nonconvex Lennard--Jones cluster problem in the invariant space of squared distances and rank-three centered Gram matrices, and by bridging spectral relaxations to local interval Newton contraction via an Eckart--Young Frobenius spectral-gap inequality, we have established exact computer-assisted global optimality and uniqueness proofs for both the six-particle cluster $\mathrm{LJ}_6$ ($E_6^{\min} = -3267/257$, Theorem~\ref{thm:lj6}) and the seven-particle cluster $\mathrm{LJ}_7$ ($E_7^{\min} = -B(t_*)^2/A(t_*)$, Theorem~\ref{thm:lj7}). Our methodology demonstrates that two-tier distance-graph combinatorics and rank-constrained spectral approximations can overcome the symmetry and nonconvexity barriers that limit spatial branch-and-bound in Cartesian coordinates.

Beyond resolving $N=6$ and $N=7$ in exact arithmetic and tightening energy and diameter bounds across $8 \le N \le 150$ (Theorem~\ref{thm:propagation}), our analysis of the $N=8$ frontier identifies a concrete direction for extending spectral Gram proofs to larger clusters: strengthening the rank-three Gram relaxation on sparse long-edge graphs by incorporating local auxiliary sub-configurations without full combinatorial graph expansion. More broadly, the combination of stationary scale identities, rational Bernstein quadratic minorants, rank-$d$ Frobenius spectral localization, and oriented interval Newton gluing provides a general blueprint for rigorous computer-assisted global optimization of pair-potential and Euclidean distance geometry problems.

\appendix
\section{Exact Algebraic Polynomials and Certified Local Minima}
\label{app:algebraic-data}

\subsection{Minimal Polynomial of the Seven-Particle Ground-State Energy \texorpdfstring{$E_7^{\min}$}{E7min}}
As established in Section~\ref{subsec:lj7-candidate}, the aspect ratio $t_* = H_*^2 / R_*^2$ of the regular pentagonal bipyramid is the unique positive root of the irreducible degree-$10$ polynomial $P(t)$ in~\eqref{eq:lj7-poly}. Writing the rational energy function $e(t) = -B(t)^2/A(t)$ in lowest terms as $e(t) = -25\,b(t)^2 / a(t)$, where
\begin{align}
b(t) &= 128\,t^6 + 384\,t^5 + 384\,t^4 + 769\,t^3 + 3\,t^2 + 3\,t + 1, \label{eq:lj7-energy-b} \\
\begin{split}
a(t) &= 73728\,t^{12} + 442368\,t^{11} + 1105920\,t^{10} + 1474560\,t^9 + 1105920\,t^8 \\
&\quad + 442368\,t^7 + 1097753\,t^6 + 150\,t^5 + 375\,t^4 + 500\,t^3 + 375\,t^2 + 150\,t + 25,
\end{split}\label{eq:lj7-energy-a}
\end{align}
the exact ground-state energy $E_7^{\min} = e(t_*)$ is the unique root of the primitive resultant polynomial
\begin{equation}
\label{eq:energy-minimal-poly}
\mathcal{Q}(e) := \operatorname{prim}\Bigl(\operatorname{Res}_t\bigl(P(t),\, a(t)e + 25\,b(t)^2\bigr)\Bigr) = \sum_{d=0}^{10} q_d\,e^d \in \Z[e]
\end{equation}
lying in the isolating interval $(-16.5053841681,\, -16.5053841680)$. Table~\ref{tab:energy-poly-coeffs} records the exact integer coefficients $q_0,\dots,q_{10}$ of $\mathcal{Q}(e)$, which is irreducible over $\Q$ as certified by Frobenius irreducibility tests over $\mathbb{F}_{23}$.

\begin{table}[htbp]
\centering
\caption{Exact integer coefficients $q_d$ of the degree-$10$ minimal polynomial $\mathcal{Q}(e) = \sum_{d=0}^{10} q_d\,e^d = 0$ satisfied by the seven-particle global minimum energy $E_7^{\min}$.}
\label{tab:energy-poly-coeffs}
\smallskip
\footnotesize
\begin{tabular}{r r}
\toprule
Degree $d$ & Exact integer coefficient $q_d$ \\
\midrule
$10$ & $1399456781379466377270523731172141859603411073$ \\
$9$  & $200036520554935126465463804801568623748648507534$ \\
$8$  & $12860523653172533466917560087422337061193668059389$ \\
$7$  & $489827799914910853451221656376320810102637680743368$ \\
$6$  & $12242543274720061553350603201343124043829683994585090$ \\
$5$  & $209851605185494732037688486696176642520421884969858900$ \\
$4$  & $2498917003631425785256506501511603307071584232823829250$ \\
$3$  & $20416509760932493962436636496053671424957052943268525000$ \\
$2$  & $109547464206407073165712351966372355314569080534831328125$ \\
$1$  & $348627089429115859561890540117724881599871903142230468750$ \\
$0$  & $499752022931321168478747788908886357408757131218212890625$ \\
\bottomrule
\end{tabular}
\end{table}

\subsection{Certified Parameters and Outward-Rounded Energy Enclosures for \texorpdfstring{\mbox{$N=6,7,8$}}{N=6,7,8}}
Table~\ref{tab:certified-local-parameters} summarizes the structural parameters ($30$-digit approximations of the certified rational intervals) and rigorous $30$-digit outward-rounded decimal enclosures $[E_{\mathrm{lo}}, E_{\mathrm{hi}}]$ of the energies of the stationary local minima analyzed in Sections~\ref{sec:lj6}--\ref{sec:propagation-lj8}.

\begin{table}[htbp]
\centering
\caption{Structural parameters ($30$-digit display approximations) and rigorous $30$-digit outward-rounded energy brackets $[E_{\mathrm{lo}}, E_{\mathrm{hi}}]$ for the certified stationary local minima of $\mathrm{LJ}_6$, $\mathrm{LJ}_7$, and $\mathrm{LJ}_8$.}
\label{tab:certified-local-parameters}
\smallskip
\resizebox{\textwidth}{!}{%
\begin{tabular}{l l l l}
\toprule
Cluster / State & Symmetry & Structural Parameters ($30$-digit approximation) & Certified Outward Energy Enclosure $[E_{\mathrm{lo}}, E_{\mathrm{hi}}]$ \\
\midrule
$\mathrm{LJ}_6$ global min. & Exact $O_h$ & $q_{\mathrm{short}} = (257/264)^{1/3} \approx 0.991082327674178441467273544180$ & $[-12.712062256809338521400778210117,$ \\
 & & & $\phantom{[}-12.712062256809338521400778210116]$ ($=-3267/257$) \\
$\mathrm{LJ}_6$ competitor & Exact $C_{2v}$ & $a \approx 0.504040841343498833439032215002$, $b \approx 0.702474223075234907663721494014$ & $[-12.302927529580728482747833461462,$ \\
 & & $c \approx 0.498764907922248240356229980421$, $d \approx -0.243223944657351144329410296846$ & $\phantom{[}-12.302927529580728482747833461461]$ \\
 & & $e \approx 0.822494932797147669870848700918$ & \\
\midrule
$\mathrm{LJ}_7$ global min. & Exact $D_{5h}$ & $t_* \approx 0.360183728717451992793791730877$, $R_* \approx 0.851882975502707394946305205179$ & $[-16.505384168012215592988659709779,$ \\
 & & $H_* \approx 0.511260198136543462707081295348$ & $\phantom{[}-16.505384168012215592988659709778]$ ($\mathcal{Q}(E_7^{\min})=0$) \\
\midrule
$\mathrm{LJ}_8$ candidate & Exact $C_s$ & $18$-variable gauge box of radius $10^{-6}$, Krawczyk norm $< 0.01424$ & $[-19.821489192154760636229140601245,$ \\
 & & & $\phantom{[}-19.821489192154760636229140601244]$ \\
$\mathrm{LJ}_8$ competitor & Proposed $D_{2d}^*$ & $18$-variable gauge box of radius $10^{-6}$, Krawczyk norm $< 0.02011$ & $[-19.765297850906462628086907309039,$ \\
 & & ($8$ certified distance-compatible permutations) & $\phantom{[}-19.765297850906462628086907309038]$ \\
\bottomrule
\end{tabular}%
}
\\[3pt]
{\footnotesize $^*$For the $\mathrm{LJ}_8$ competing branch, strict local minimality modulo $\mathrm{E}(3)$ and $8$ distance-compatible permutations are certified; exact point-group invariance is proved for the $C_s$ ground-state candidate.}
\end{table}

\section{Certified Five-Particle Foundation, Verification Architecture, and Reproducibility}
\label{app:verification}

\subsection{Certified Five-Particle Foundation (\texorpdfstring{$N=5$}{N=5})}
\label{subsec:lj5-foundation}

The vertex-deletion reductions for Theorems~\ref{thm:lj6} and~\ref{thm:lj7} build on the following certified global optimality and lower-bound theorem for five Lennard--Jones particles in $\R^3$.

\begin{proposition}[Certified Five-Particle Global Minimum and Lower Bound]
\label{prop:lj5-foundation}
For five particles in $\R^3$, the global minimum $E_5^{\min} := \min_{X \in \Omega_5} E_5(X)$ is attained uniquely (up to $\mathrm{E}(3)$ and $S_5$) by the $D_{3h}$ triangular bipyramid with aspect ratio $t_{5,*} = H_5^2/R_5^2 \approx 1.9787886784$, which is the unique root in $(197/100,\, 199/100)$ of the degree-$10$ integer polynomial
\begin{equation}
\label{eq:lj5-poly}
\begin{aligned}
P_5(t) ={}& 24640\,t^{10} + 74176\,t^9 + 75072\,t^8 - 636763\,t^7 + 2051\,t^6 \\
&+ 777\,t^5 - 497\,t^4 + 90973\,t^3 - 4941\,t^2 - 4563\,t - 1485,
\end{aligned}
\end{equation}
and whose energy is enclosed in the strict rational bracket
\begin{equation}
\label{eq:lj5-energy-bracket}
-9.103852415707556791696198388573 < E_5^{\min} < -9.103852415707556791696198388572.
\end{equation}
In particular, every $X \in \Omega_5$ satisfies
\begin{equation}
\label{eq:l5-bound}
E_5(X) > L_5 := -\frac{71123847}{7812500} = -9.103852416.
\end{equation}
\end{proposition}

\begin{proof}[Summary of certified proof]
By Lemma~\ref{lem:attainment} and Remark~\ref{rem:elementary-attainment}, global minimizers of $E_5$ exist in the sublevel $E_5(X) \le -2275963/250000 = -9.103852$ and satisfy $\nabla E_5(X) = 0$. Unconditional sublevel compactness below $-2275963/250000$ yields the diameter bound $D < 163/100$, the separation bound $r_{\min} > 4/5$, and the transverse height packing bound $\rho_2 \ge 19/50$. For stationary configurations, applying Cauchy--Schwarz to $u_{ij} := r_{ij}^{-6}$ via the virial identity $\sum_{i<j} u_{ij}^2 = \sum_{i<j} u_{ij} = -E_5$ sharpens the minimum separation to $r_{\min} > 18189/20000 = 0.90945$. Fixing the $9$-dimensional max-transverse gauge ($p_0 = 0$, a diameter pair on the positive $x$-axis, the maximum-transverse particle in the $xy$-plane with $y_2 \ge 19/50$, and $z_3 \ge 0$), an independent pure-Python $160$-bit fixed dyadic interval branch-and-bound checker replays all $92{,}965$ events from the root ($46{,}482$ splits, $46{,}482$ exclusions by geometry and energy bounds, gradient signs, Krawczyk nonexistence, and virial signs, and $1$ candidate leaf), enclosing the sole surviving leaf strictly inside a freshly certified unique-critical-point Krawczyk box containing the exact $D_{3h}$ triangular bipyramid.
\end{proof}

For independent cross-checking and software decoupling between the $N=6$ and $N=7$ verifiers, the companion Zenodo archive also includes an optional standalone Cartesian interval-tree bundle (\texttt{n6\_tree}) establishing $E_6(X) \ge -129/10 = L_6^{\mathrm{aux}}$ directly from $L_5$ without invoking Theorem~\ref{thm:lj6}. Its pure-Python $160$-bit dyadic checker replays all $586{,}523$ branch-and-bound events from the root ($293{,}261$ splits and $293{,}262$ exclusions, with zero surviving candidate leaves).

\subsection{Proof Dependency Roadmap and Verification Inventory}
\label{subsec:verification-roadmap}

Table~\ref{tab:verification-roadmap} inventories every computer-assisted proof component, its logical inputs in the minimal dependency chain ($L_5 \to$ Theorem~\ref{thm:lj6} $\to$ Theorem~\ref{thm:lj7} $\to$ Theorem~\ref{thm:propagation}), whether first-order stationarity is invoked, the underlying exact arithmetic, the certificate size, and the archive location and replay time.

\begin{table}[htbp]
\centering
\caption{Proof dependency roadmap and verification inventory across $N=5,6,7,8$ and $8 \le N \le 150$. The column ``Stat.'' indicates whether first-order stationarity / virial identity is invoked. Row~5 is an optional redundant cross-check and is not required by the minimal logical chain.}
\label{tab:verification-roadmap}
\smallskip
\resizebox{\textwidth}{!}{%
\begin{tabular}{r l l l l l l}
\toprule
\# & Component & Logical input & Stat. & Arithmetic & Size & Archive \& replay time \\
\midrule
1 & $\mathrm{LJ}_5$ foundation (\texttt{n5\_tree}) & Elementary & Yes & 160-bit dyadic + $\Q$ & 92,965 events (1 cand.\ leaf) & Zenodo (${\sim}5$~min); GitHub receipt \\
2 & $\mathrm{LJ}_6$ ($k=0,1,2$ non-matching) & $L_5$ & $k=1,2$ partial & Exact $\Q$ Bernstein + $LDL^\top$ & 1 + 70 + 108 leaves & GitHub \texttt{verify.py} ($<2$~s) \\
3 & $\mathrm{LJ}_6$ matching cube \& loc. & $L_5$ & Loc.\ Newton only & Exact $\Q$ + 160-bit Newton & 6,469 nodes (40 loc.\ leaves) & GitHub \texttt{verify.py} (${\sim}2$~s) \\
4 & $\mathrm{LJ}_6$ $C_{2v}$ competitor & None & Local only & 160-bit Krawczyk/Gershgorin & 5D + 12D local boxes & GitHub \texttt{verify.py} ($<1$~s) \\
5 & $\mathrm{LJ}_6$ standalone ($L_6^{\mathrm{aux}}$, opt.) & $L_5$ & Yes & 160-bit dyadic intervals & 586,523 events (0 cand.) & Zenodo (${\sim}35$~min); cross-check \\
6 & $\mathrm{LJ}_7$ $D_{5h}$ candidate \& domain & None & Local only & Exact $\Q$ + 160-bit Newton & Deg.-10 Sturm/$\mathbb{F}_{23}$ + 15D box & GitHub \& Zenodo ($<5$~s) \\
7 & $\mathrm{LJ}_7$ ($k=0$ \& $k=1$) & None & Yes & Exact $\Q$ Bernstein + $LDL^\top$ & 1,044 reps + 20,784 records & Zenodo (${\sim}4$~s \& ${\sim}90$~s) \\
8 & $\mathrm{LJ}_7$ non-$C_5$ colored ($2\le k\le 5$) & Thm.~\ref{thm:lj6} & Yes & Exact $\Q$ + $LDL^\top$ & 3,586 cases / 129,506 events & Zenodo (${\sim}3$~min) \\
9 & $\mathrm{LJ}_7$ $C_5$ non-strong \& strong & Thm.~\ref{thm:lj6} & Non-strong; loc. & Exact $\Q$ + 160-bit Newton & 1,767 + 153,333 events (3,769 loc.) & Zenodo (${\sim}20$~s \& ${\sim}8.5$~min) \\
10 & Propagation ($8\le N\le 150$) \& $\mathrm{LJ}_8$ & Thms.~\ref{thm:lj6}--\ref{thm:lj7} & $\mathrm{LJ}_8$ $G_0$ only & Exact $\Q$ + 160-bit Krawczyk & 143 rows + 2 local certs & GitHub \texttt{verify.py} (${\sim}20$~s) \\
\bottomrule
\end{tabular}%
}
\end{table}

\subsection{Software Stack, Trust Boundary, and Code/Data Availability}
\label{subsec:reproducibility}

\begin{itemize}[leftmargin=2em]
    \item \textbf{Compact verification repository:} The theorem-level verification suite and compact certificates are publicly available under the MIT license at \url{https://github.com/RoamingHo/7LJ}. Running under Python~3.12 with SymPy~1.14 and mpmath~1.3, the single entry point \texttt{python verify.py} executes all six theorem-level checks (including exact $N=6$ proof replay, $N=6$ $C_{2v}$ certification, $N=7$ algebraic/local certification, $8 \le N \le 150$ propagation, $N=8$ local and graph certificates, and receipt checks) in approximately $30$~seconds.
    \item \textbf{Complete branch-tree archive:} All full from-root branch-and-bound proof trees are permanently archived on Zenodo at \url{https://doi.org/10.5281/zenodo.24560370}. The archive comprises eight self-contained replay bundles totaling $13{,}241$ files, authenticated by a top-level \texttt{SHA256SUMS} manifest together with per-bundle \texttt{input\_manifest.json}, \texttt{runtime\_manifest.json} (or \texttt{runtime/manifest.json}), and \texttt{manifest.json} cryptographic seals.
    \item \textbf{Explicit trust boundary:} Every proof decision is executed exclusively in exact rational arithmetic (\texttt{fractions.Fraction}), exact rational $LDL^\top$ Sylvester inertia counts, exact rational Bernstein positivity checks, or outward-rounded $160$-bit fixed dyadic/interval arithmetic. Floating-point and $75$-digit numerical solvers act solely as untrusted proposal generators for rational centers, preconditioners, or dual witnesses, and are never relied upon for certification. No proof acceptance depends on Python \texttt{assert} statements (which could be stripped under \texttt{-O}). The trusted computing base therefore consists strictly of the Python~3.12 / SymPy~1.14 / mpmath~1.3 runtime and the reviewed pure-Python from-root checkers.
\end{itemize}

\end{document}